\documentclass{amsart}
\usepackage{amsmath}
\usepackage{mathrsfs}
\usepackage{amsfonts}
\usepackage{amssymb}
\usepackage{graphicx}
\usepackage{amsthm}
\usepackage[toc,page]{appendix}
\usepackage{geometry}
\usepackage[all]{xy}
\usepackage{lmodern}
\usepackage[T1]{fontenc}
\usepackage[dvipsnames]{xcolor}
\usepackage{graphics}
\usepackage{comment}
\usepackage{csquotes}
\usepackage[colorlinks=true,linkcolor=NavyBlue,urlcolor=RoyalBlue,citecolor=PineGreen,linktocpage=true]{hyperref}
\usepackage{enumerate}
\usepackage{caption}
\usepackage{subcaption}
\usepackage{float}
\usepackage{multicol}

\numberwithin{equation}{section}

\newtheorem{Theorem}{Theorem}[section]	
\newtheorem*{Theorem*}{Theorem}
\newtheorem{proposition}[Theorem]{Proposition} 
\newtheorem{lemma}[Theorem]{Lemma}
	
\newtheorem{defn}[Theorem]{Definition}
\newtheorem*{defn*}{Definition}
\newtheorem{Cor}[Theorem]{Corollary}
\newtheorem*{remark}{Remark}
\newtheorem*{remarks}{Remarks}
\newtheorem{example}{Example}[section]

\newcommand{\supp}{{\operatorname{supp}}}

\DeclareTextFontCommand{\Emph}{\bfseries\itshape}

\usepackage{tikz}
\usetikzlibrary{decorations,arrows}
\usetikzlibrary{decorations.pathmorphing}
\usepgflibrary{decorations.pathreplacing} 
\usepackage{pgf}
\usetikzlibrary{patterns}

\tikzset{
	schraffiert/.style={pattern=horizontal lines,pattern color=#1},
	schraffiert/.default=black
}

\tikzset{
	ultra thin/.style= {line width=0.1pt},
	very thin/.style=  {line width=0.2pt},
	thin/.style=       {line width=0.4pt},
	semithick/.style=  {line width=0.6pt},
	thick/.style=      {line width=0.8pt},
	very thick/.style= {line width=1.2pt},
	ultra thick/.style={line width=2.4pt}
}

\usetikzlibrary{shapes}
\usetikzlibrary{positioning}

\title{Automaticity of structural defects in fixed points of block substitutions}
\author{Lior Tenenbaum}
\date{}

\address{Department of Mathematics\\
	Bar-Ilan university\\
	Ramat Gan, Israel}
\email{lior.tenenbaum@biu.ac.il}

\begin{document}
\begin{abstract}
		We study structural defects in fixed points of block substitutions, building on the work of Band, Beckus, Pogorzelski, and Tenenbaum (2026)
on spectral pollution in substitution systems, and providing its
technical supplement. We develop an automata-based framework to characterize their occurrence and derive verifiable combinatorial criteria. We also introduce simplifications for  the construction and analysis of these automata, including tensor-product decompositions, that  clarify the mechanisms underlying structural defects and their emergence.
	\end{abstract}
\maketitle


\section{Introduction}
Recently, Beckus, Bellissard, and coauthors used approximations of the underlying dynamical structures to study spectra of more complicated operators; see
\cite{BeckusThesis,BBdN18,BBC19,BBdN20,BecTak21,BBT24}. 
In a sequence of works
\cite{BeckusThesis,BBdN20,Ten24,BBPT24,BBPT26},
this approach was applied to spectra arising from substitution systems. 
More specifically, starting from an initial periodic configuration $\omega_0$, one iterates the substitution on $\omega_0$ and considers the dynamical system obtained as the closure of its translation orbit. The evolution of these iterated orbit closures was studied in \cite{BBPT24,BBPT26}, while a description of the corresponding limiting dynamical systems was obtained in \cite{BBPT26}.
Following \cite{BBPT26}, where limiting structural defects in substitution systems and their spectral implications were studied, we further investigate the combinatorial conditions governing the appearance of structural defects. While these conditions were implicit in \cite{BBPT26}, we make them explicit here. Accordingly, this note may be viewed as a technical supplement to \cite{BBPT26}, clarifying the combinatorial assumptions underlying the structural-defect criteria of \cite[Section~2]{BBPT26}.
\newline
A patch $P$ is called $S$-legal, with respect to a block substitution $S$, if there exists some letter $a\in \mathcal{A}$ and $n \in \mathbb{N}$ such that $P$ is a subpatch of $S^n(a)$. A patch $P$ occurring in a limit point configuration $\rho= \lim_{n\to \infty} S^n(\omega)$, for some starting configuration $\omega$, is called a structural defect (with respect to $S$) if it is not $S$-legal. When a configuration $\rho \in \mathcal{A}^{\mathbb{Z}^d}$ satisfies $S^k(\rho)=\rho$ for some $k\in \mathbb{N}$, we call it a fixed point of $S$. Restricting attention to the set $T=\{-1,0\}^d$, every fixed point $\rho$ is uniquely determined by its restriction $\rho|_{T}$ together with the substitution $S$, as was described in \cite[Section 4]{BBPT26}.

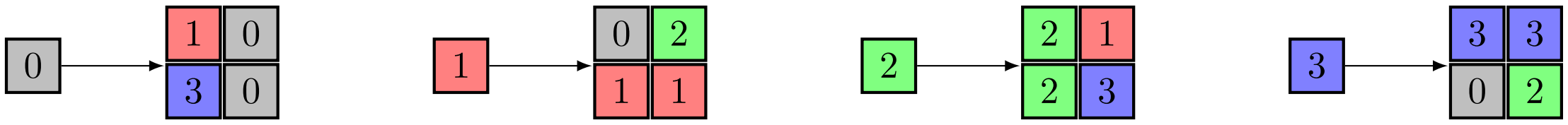
\begin{figure}[H]
			\begin{tikzpicture}
				[
				node distance= 0.1mm, 
				0node/.style={rectangle, draw=black, fill=gray!50!white, thick, minimum size=5mm, font ={$0$} },
				1node/.style={rectangle, draw=black, fill=red!50!white, thick, minimum size=5mm, font ={$1$}},
				2node/.style={rectangle, draw=black, fill=green!50!white, thick, minimum size=5mm, font ={$2$}},
				3node/.style={rectangle, draw=black, fill=blue!50!white, thick, minimum size=5mm, font={$3$}},
				]

				\begin{scope}
					
					\node[0node] (0Loc) [yshift=-0.3cm] {};
					
					\begin{scope}[xshift=1.5cm ,yshift=0cm]
						\node  (0Im) [yshift=-0.3cm ,xshift=-0.15cm] {};
						\node[1node] (LU) {};
						\node[0node] (RU) [right= of LU] {};
						\node[3node] (LB) [below= of LU] {};
						\node[0node] [right= of LB] {};
					\end{scope}
					
				\end{scope}

				\begin{scope}[xshift= 4cm]
					
					\node[1node] (1Loc) [yshift=-0.3cm] {};
					
					\begin{scope}[xshift=1.5cm ,yshift=0cm]
						\node  (1Im) [yshift=-0.3cm ,xshift=-0.15cm] {};
						\node[0node] (LU) {};
						\node[2node] (RU) [right= of LU] {};
						\node[1node] (LB) [below= of LU] {};
						\node[1node] [right= of LB] {};
					\end{scope}
					
				\end{scope}
				
				\begin{scope}[xshift= 8cm]
					
					\node[2node] (2Loc) [yshift=-0.3cm] {};
					
					\begin{scope}[xshift=1.5cm ,yshift=0cm]
						\node  (2Im) [yshift=-0.3cm ,xshift=-0.15cm] {};
						\node[2node] (LU) {};
						\node[1node] (RU) [right= of LU] {};
						\node[2node] (LB) [below= of LU] {};
						\node[3node] [right= of LB] {};
					\end{scope}
					
				\end{scope}

				\begin{scope}[xshift=12cm]
					
					\node[3node] (3Loc) [yshift=-0.3cm] {};
					
					\begin{scope}[xshift=1.5cm ,yshift=0cm]
						\node  (3Im) [yshift=-0.3cm ,xshift=-0.15cm] {};
						\node[3node] (LU) {};
						\node[3node] (RU) [right= of LU] {};
						\node[0node] (LB) [below= of LU] {};
						\node[2node] [right= of LB] {};
					\end{scope} 
					
				\end{scope}

				\draw[-latex]  (0Loc) edge (0Im) ;
				\draw[-latex]  (1Loc) edge (1Im) ;
				\draw[-latex]  (2Loc) edge (2Im) ;
				\draw[-latex]  (3Loc) edge (3Im) ;

			\end{tikzpicture}
			\caption{The substitution rule for the symbolic table substitution.}	 
			\label{Fig:TableRule}	 
		\end{figure}

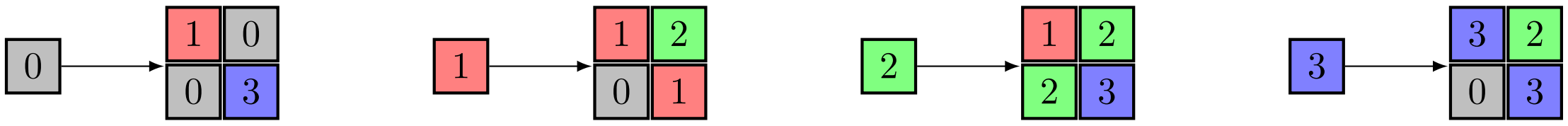
\begin{figure}[H]
			\begin{tikzpicture}
				[
				node distance= 0.1mm, 
				0node/.style={rectangle, draw=black, fill=gray!50!white, thick, minimum size=5mm, font ={$0$} },
				1node/.style={rectangle, draw=black, fill=red!50!white, thick, minimum size=5mm, font ={$1$}},
				2node/.style={rectangle, draw=black, fill=green!50!white, thick, minimum size=5mm, font ={$2$}},
				3node/.style={rectangle, draw=black, fill=blue!50!white, thick, minimum size=5mm, font={$3$}},
				]

				\begin{scope}
					
					\node[0node] (0Loc) [yshift=-0.3cm] {};
					
					\begin{scope}[xshift=1.5cm ,yshift=0cm]
						\node  (0Im) [yshift=-0.3cm ,xshift=-0.15cm] {};
						\node[1node] (LU) {};
						\node[0node] (RU) [right= of LU] {};
						\node[0node] (LB) [below= of LU] {};
						\node[3node] [right= of LB] {};
					\end{scope}
					
				\end{scope}

				\begin{scope}[xshift= 4cm]
					
					\node[1node] (1Loc) [yshift=-0.3cm] {};
					
					\begin{scope}[xshift=1.5cm ,yshift=0cm]
						\node  (1Im) [yshift=-0.3cm ,xshift=-0.15cm] {};
						\node[1node] (LU) {};
						\node[2node] (RU) [right= of LU] {};
						\node[0node] (LB) [below= of LU] {};
						\node[1node] [right= of LB] {};
					\end{scope}
					
				\end{scope}
				
				\begin{scope}[xshift= 8cm]
					
					\node[2node] (2Loc) [yshift=-0.3cm] {};
					
					\begin{scope}[xshift=1.5cm ,yshift=0cm]
						\node  (2Im) [yshift=-0.3cm ,xshift=-0.15cm] {};
						\node[1node] (LU) {};
						\node[2node] (RU) [right= of LU] {};
						\node[2node] (LB) [below= of LU] {};
						\node[3node] [right= of LB] {};
					\end{scope}
					
				\end{scope}

				\begin{scope}[xshift=12cm]
					
					\node[3node] (3Loc) [yshift=-0.3cm] {};
					
					\begin{scope}[xshift=1.5cm ,yshift=0cm]
						\node  (3Im) [yshift=-0.3cm ,xshift=-0.15cm] {};
						\node[3node] (LU) {};
						\node[2node] (RU) [right= of LU] {};
						\node[0node] (LB) [below= of LU] {};
						\node[3node] [right= of LB] {};
					\end{scope} 
					
				\end{scope}

				\draw[-latex]  (0Loc) edge (0Im) ;
				\draw[-latex]  (1Loc) edge (1Im) ;
				\draw[-latex]  (2Loc) edge (2Im) ;
				\draw[-latex]  (3Loc) edge (3Im) ;

			\end{tikzpicture}
			\caption{The substitution rule for the symbolic chair substitution.}	 
			\label{Fig:ChairRule}	 
		\end{figure}

In \cite{BBPT26},  the existence of structural defects in fixed points associated with the symbolic table and symbolic chair tilings was studied using ad hoc methods. We briefly recall these two block substitutions, originally introduced in \cite{Rob99}. Both are two-dimensional $2\times 2$ substitutions over the four letter alphabet $\mathcal{A}=\{0,1,2,3\}$.  The substitution rule for the symbolic table tiling substitution $S_{\mathtt{table}}$  is shown in Figure \ref{Fig:TableRule}.
	The rule for the symbolic chair tiling substitution $S_{\mathtt{Chair}}$  is given in Figure \ref{Fig:ChairRule}.
	In \cite[Corollary 2.10]{BBPT26} it was shown that constant letter configurations generate pure-line structural defects. That is, every $2\times 2$ patch occurring  in the set 
$$ \Lambda_1:=\{ (m_1,m_2): m_2\in \{-1,0\}, m_1\in \mathbb{Z} \}= \mathbb{Z}\mathsf{e}_1+T $$
or
$$\Lambda_2=\{ (m_1,m_2): m_1\in \{-1,0\}, m_2\in \mathbb{Z} \} = \mathbb{Z}\mathsf{e}_2+T $$
is illegal. 
The significance of these pure-line defects for spectral approximation was demonstrated in \cite[Proposition 2.8 and Theorem 2.9]{BBPT26}, where it was shown that they may lead spectral pollution of positive Lebesgue measure. However, explicit criteria guaranteeing their occurrence were not provided there. To study their appearance, we show that they are governed by automata defined on either $2\times 1$-patches or $1\times 2$-patches. 
The construction of these automata is described in Section \ref{Sect:NumCode} and the mechanism underlying these automata is discussed in Section ~\ref{Sect:Defect}.

\begin{figure}[b]
	\centering
	\begin{tikzpicture}
				[
				node distance= 0.1mm, 
				0node/.style={rectangle, draw=black, fill=gray!50!white, thick, minimum size=5mm, font ={$0$} },
				1node/.style={rectangle, draw=black, fill=red!50!white, thick, minimum size=5mm, font ={$1$}},
				2node/.style={rectangle, draw=black, fill=green!50!white, thick, minimum size=5mm, font ={$2$}},
				3node/.style={rectangle, draw=black, fill=blue!50!white, thick, minimum size=5mm, font={$3$}},
				]
				\node[1node] (0U) {};
				\node[1node] (0B) [below= of 0U]  {};
				\node[1node] (-1U) [left= of 0U]  {};
				\node[1node] (-1B) [below= of -1U]  {};
				
				\node[1node] (1U) [right= of 0U]  {};
				\node[0node] (1B) [below= of 1U]  {};
				\node[1node] (2U) [right= of 1U]  {};
				\node[2node] (2B) [below= of 2U]  {};
				\node[1node] (3U) [right= of 2U]  {};
				\node[1node] (3B) [below= of 3U]  {};
				
				\node[1node] (4U) [right= of 3U]  {};
				\node[0node] (4B) [below= of 4U]  {};
				\node[1node] (5U) [right= of 4U]  {};
				\node[2node] (5B) [below= of 5U]  {};
				\node[1node] (6U) [right= of 5U]  {};
				\node[1node] (6B) [below= of 6U]  {};
				
				\node[1node] (7U) [right= of 6U]  {};
				\node[0node] (7B) [below= of 7U]  {};
				\node[1node] (8U) [right= of 7U]  {};
				\node[2node] (8B) [below= of 8U]  {};
				\node[1node] (9U) [right= of 8U]  {};
				\node[1node] (9B) [below= of 9U]  {};
				
				\node[1node] (10U) [right= of 9U]  {};
				\node[0node] (10B) [below= of 10U]  {};
				\node[1node] (11U) [right= of 10U]  {};
				\node[2node] (11B) [below= of 11U]  {};
				\node[1node] (12U) [right= of 11U]  {};
				\node[1node] (12B) [below= of 12U]  {};	
				
				\node[1node] (13U) [right= of 12U]  {};
				\node[0node] (13B) [below= of 13U]  {};
				
				\node[1node] (-2U) [left= of -1U]  {};
				\node[0node] (-2B) [below= of -2U]  {};
				\node[1node] (-3U) [left= of -2U]  {};
				\node[2node] (-3B) [below= of -3U]  {};
				\node[1node] (-4U) [left= of -3U]  {};
				\node[1node] (-4B) [below= of -4U]  {};	
				
				\node[1node] (-5U) [left= of -4U]  {};
				\node[0node] (-5B) [below= of -5U]  {};
				\node[1node] (-6U) [left= of -5U]  {};
				\node[2node] (-6B) [below= of -6U]  {};
				\node[1node] (-7U) [left= of -6U]  {};
				\node[1node] (-7B) [below= of -7U]  {};	
				
				\node[1node] (-8U) [left= of -7U]  {};
				\node[0node] (-8B) [below= of -8U]  {};
				\node[1node] (-9U) [left= of -8U]  {};
				\node[2node] (-9B) [below= of -9U]  {};
				\node[1node] (-10U) [left= of -9U]  {};
				\node[1node] (-10B) [below= of -10U]  {};	
				
				\node[1node] (-11U) [left= of -10U]  {};
				\node[0node] (-11B) [below= of -11U]  {};
				\node[1node] (-12U) [left= of -11U]  {};
				\node[2node] (-12B) [below= of -12U]  {};
				\node[1node] (-13U) [left= of -12U]  {};
				\node[1node] (-13B) [below= of -13U]  {};
				
				\node[1node] (-14U) [left= of -13U]  {};
				\node[0node] (-14B) [below= of -14U]  {};
				
	\end{tikzpicture}
	\caption{\centering A portion of a pure-line defect supported on $\Lambda_1$ in a fixed point of the symbolic table tiling substitution.}
	\label{Fig:PureLineDefectConst1}
\end{figure}
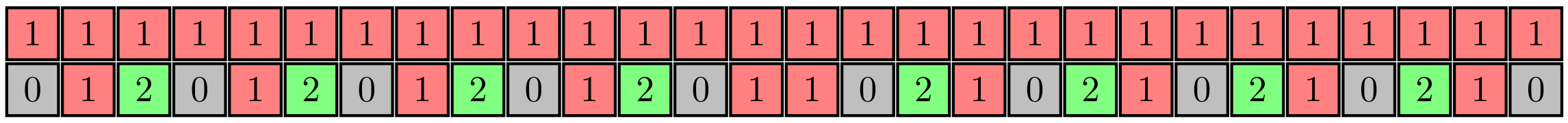

This automaton structure is natural in light of Cobham's theorem, see \cite{Cob72} and \cite[Chapter 6.3]{Shal03}. Roughly speaking, the theorem states that automatic sequences are equivalent to fixed points of constant-length substitutions. The main distinction in our setting is that we consider fixed points that are not necessarily canonically associated with the substitution, and our automata describe the evolution of larger patches rather than individual symbols.
\newline 
Using these automata, we obtain in Lemma~\ref{Lem:PL-Defect} a necessary and sufficient criterion for the existence of pure-line defects. Furthermore, Proposition~\ref{Prop:PureLineDefect} provides an easily verifiable sufficient condition for their emergence, generalizing the behavior observed for the symbolic table tiling substitution.
\newline   
On the other hand, \cite[Corollary 2.3]{BBPT26} shows that every iterated substitutive approximation of the symbolic chair tiling admits a subsequence converging to the target spectrum, corresponding to the ideal substitution system, up to a set of measure zero and without essential spectral pollution. More generally, \cite[Theorem 2.2 and Lemma 5.9]{BBPT26} proves that the same behavior holds for purely $T$-illegal substitutions, introduced in \cite[Definition 5.7]{BBPT26}. Although it was shown in \cite[Lemma A.1]{BBPT26} that the symbolic chair tiling substitution is purely $T$-illegal, Corollary~\ref{Cor:BddDefect} provides a more general criterion for verifying this property using the automaton description of structural defects in fixed points. As shown in \cite{BBPT26}, this behavior occurs for all one-dimensional symbolic substitution systems.
\newline 
To facilitate the analysis of the automata introduced in this note, Section~\ref{Sect:Tensor} presents a tensor-product decomposition that breaks them into smaller components and thereby reduces the computational complexity of the relevant computations.

We note that there are several natural avenues to generalize the results in this note. These include computing structural  defects for quasi-fixed points and for quasi-fixed points arising in tile-substitutions; see \cite{Kraw26}. For the sake of simplicity we only discuss structural defects on fixed points of symbolic block substitutions.
\newline
The note is organized as follows. In Section~\ref{Sect:Prelim}, we introduce the necessary preliminaries, including block substitutions and their associated automata. In Section~\ref{Sect:Defect}, we discuss combinatorial conditions on substitutions that govern the occurrence of various types of structural defects, including those required in \cite[Section~2]{BBPT26}. Finally in Section \ref{Sect:Simpl}, we present graph-based simplifications of the automata that reduce the computational effort required to construct them.

\subsection*{AI Disclosure}

The author used AI language models, including ChatGPT and Claude, to assist with proofreading, manuscript typesetting, and phrasing suggestions. The author takes full responsibility for the final text.

\section{Preliminaries} \label{Sect:Prelim}

Let $\mathcal{A}$ be a finite nonempty set, called an \Emph{alphabet}, and consider the ambient space $\mathbb{Z}^d$ for some $d\in \mathbb{N}$. Write
\[ \operatorname{Pat}(\mathcal{A}^{\mathbb{Z}^d})=\{ P:M\to \mathcal{A}: \; M\subseteq \mathbb{Z}^d \} \]
for the set of patches. For a patch $P \in \operatorname{Pat}(\mathcal{A}^{\mathbb{Z}^d})$, $\supp(P):=M$ denotes its support. When $M=\mathbb{Z}^d$ it is called a configuration, and it is also sometimes called a partial configuration when $M\neq \mathbb{Z}^d$.  There is a natural $\mathbb{Z}^d$-action induced on $\operatorname{Pat}(\mathcal{A}^{\mathbb{Z}^d})$ by
\[ [\mathsf{m}+P](\mathsf{k}):= P(\mathsf{k}-\mathsf{m}) \quad \text{for} \quad \mathsf{m}, \mathsf{k}\in \mathbb{Z}^d. \]
For two patches $P,Q\in \operatorname{Pat}(\mathcal{A}^{\mathbb{Z}^d})$, we say $Q$ is a \Emph{subpatch} of $P$, and write $Q\prec P$, if
\[ \mathsf{m} + \supp(P)\subseteq \supp(Q) \quad \text{and} \quad \mathsf{m}+P= Q|_{\mathsf{m}+\supp(P)}. \]
For a configuration $\omega : \mathbb{Z}^d \to \mathcal{A}$, its set of associated patches is defined by
\[ W(\omega):= \{ P\in \operatorname{Pat} (\mathcal{A}^{\mathbb{Z}^d}): P = \omega |_{M}, \; M\subseteq \mathbb{Z}^d\}. \]
 We call this set of patches is called the \Emph{dictionary} of the configuration $\omega$.

We consider the set $T=\{-1,0\}^d$, which was extensively studied in \cite{BBPT24}. For every $1\leq j \leq d$, we define the set
\begin{equation} \label{Eq:TjDefn}
T^j:= T \cap \{  (x_1,...,x_d)\in \mathbb{R}^d: x_j=0 \}.
\end{equation}
We denote the standard basis vectors of $\mathbb{R}^d$ by $\mathsf{e}_1,...,\mathsf{e}_d$, and the set of nonnegative integers by $\mathbb{N}_0:=\mathbb{N}\cup\{0\}$. We fix a norm $|| \cdot ||$ on $\mathbb{R}^d$, and denote for $r>0$ the ball $B_r:=\{ x\in \mathbb{R}^d: ||x||<r  \}$.

\subsection{Symbolic block substitutions}
We consider in this note only block substitutions over $\mathbb{Z}^d$, see \cite{Fogg02, NPF05, BaakeGrimm13, Bart18}. Given a finite alphabet $\mathcal{A}$, a substitution of length $\mathsf{m}\in \mathbb{N}^d$ is determined by a substitution rule $S_0:\mathcal{A}\to K_{\mathsf{m}}$ which extends locally to finite patches and configurations. 
Our description of block substitutions follows the framework of \cite{NPF05}.
Given a vector $\mathsf{q}=(q_1,...,q_d)\in \mathbb{N}^d$, we consider the dilation mapping
\[ D_{\mathsf{q}}:\mathbb{R}^d\to \mathbb{R}^d, \quad  D_{\mathsf{q}}(x_1,...,x_d)=(q_1x_1,...,q_dx_d). \]
We note that $D_{\mathsf{q}}[\mathbb{Z}^d]\subseteq \mathbb{Z}^d$. Given  $\mathsf{q}\in \mathbb{N}^d$ and a finite alphabet $\mathcal{A}$, a substitution rule of length $\mathsf{q}$ over $\mathcal{A}$ is a mapping $S_0: \mathcal{A}\to \mathcal{A}^{K_{\mathsf{q}}}$ where $K_{\mathsf{q}}:=\prod_{j=1}^d \{0,..., q_j-1\}$. The fundamental cell $[0,1)^d$ satisfies
\[ D_{\mathsf{q}}\big[ [0,1)^d \big] \cap \mathbb{Z}^d = K_{\mathsf{q}} \quad \text{and} \quad  \bigsqcup_{ \mathsf{m}\in \mathbb{Z}^d } \big( \mathsf{m}+ [0,1)^d \big) = \mathbb{R}^d.  \]
Moreover, $D_\mathsf{q}$ is an automorphism of $\mathbb{R}^d$ as an additive group. 
The substitution rule is extended to an arbitrary patch $P:F\to \mathcal{A}$ via
\[ [S(P)]\big( D_\mathsf{q}(\mathsf{m})+\mathsf{r}\big):=  \big[ S_0\big( P(\mathsf{m}) \big) \big](\mathsf{r}) \quad \text{for} \quad \mathsf{m}\in F,\; \mathsf{r}\in K_{\mathsf{q}}, \]
where $\operatorname{supp}\big( S(P) \big)=  D_\mathsf{q}\Big[ F+[0,1)^d \Big]\cap \mathbb{Z}^d$. A simple induction argument shows that
\begin{equation} \label{Eq:SuppIter}
\operatorname{supp}\big( S^n(P) \big)= D_\mathsf{q}^n\big[ \operatorname{supp}(P)+[0,1)^d  \big] \cap \mathbb{Z}^d \quad \text{for any} \;\; n\in \mathbb{N}. 
\end{equation}
We implicitly assume throughout that the vector $\mathsf{q}$ satisfies that $q_j\geq 2$ for all $1\leq j \leq d $. Under this assumption, $\lim_{n\to \infty} \operatorname{supp}\big( S^n(P) \big)=\mathbb{Z}^d$, whenever $T\subseteq\supp(P)$.
Let us consider the following example of a block substitution.
\begin{example}
Let $\mathcal{A}=\{0,1\}$ be the binary alphabet, let $D_3$ be the square dilation map $D_3(x_1,x_2)=(3x_1,3x_2)$, and let the symbolic squiral substitution rule be as given in Figure \ref{Fig:Squiral_Rule}. 
\begin{figure}[ht]\centering 
			\begin{tikzpicture}
				[
				node distance= 0.1mm, 
				0node/.style={rectangle, draw=black, fill=gray!50!white, thick, minimum size=5mm, font ={$0$} },
				1node/.style={rectangle, draw=black, fill=red!50!white, thick, minimum size=5mm, font ={$1$}},
				2node/.style={rectangle, draw=black, fill=green!50!white, thick, minimum size=5mm, font ={$2$}},
				3node/.style={rectangle, draw=black, fill=blue!50!white, thick, minimum size=5mm, font={$3$}},
				]

				\begin{scope}
					
					\node[0node] (0Loc) [yshift=-0.5cm] {};
					
					\begin{scope}[xshift=1.5cm ,yshift=0cm]
						\node  (0Im) [yshift=-0.5cm ,xshift=-0.15cm] {};
						\node[1node] (LU) {};
						\node[0node] (MU) [right= of LU] {};
						\node[1node] (RU) [right= of MU] {};
						\node[0node] (LM) [below= of LU] {};
						\node[0node] (MM) [right= of LM] {};
						\node[0node] (RM) [right= of MM] {};
						\node[1node] (LB) [below= of LM] {};
						\node[0node] (MB) [right= of LB] {};
						\node[1node] [right= of MB] {};
					\end{scope}
					
				\end{scope}

				\begin{scope}[xshift= 5cm]
					
					\node[1node] (1Loc) [yshift=-0.5cm] {};
					
					\begin{scope}[xshift=1.5cm ,yshift=0cm]
						\node  (1Im) [yshift=-0.5cm ,xshift=-0.15cm] {};
						\node[0node] (LU) {};
						\node[1node] (MU) [right= of LU] {};
						\node[0node] (RU) [right= of MU] {};
						\node[1node] (LM) [below= of LU] {};
						\node[1node] (MM) [right= of LM] {};
						\node[1node] (RM) [right= of MM] {};
						\node[0node] (LB) [below= of LM] {};
						\node[1node] (MB) [right= of LB] {};
						\node[0node] [right= of MB] {};
					\end{scope}
					
				\end{scope}

				\draw[-latex]  (0Loc) edge (0Im) ;
				\draw[-latex]  (1Loc) edge (1Im) ;

			\end{tikzpicture}
			\caption{The substitution rule for the symbolic squiral substitution.}	 
			\label{Fig:Squiral_Rule}
\end{figure}
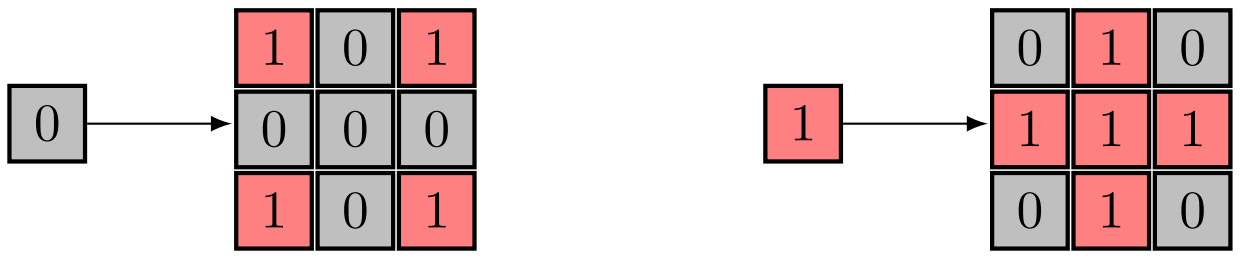
The symbolic squiral substitution was previously studied in \cite{BG14}. Note that $S_{\mathtt{Squir}}$ is bijective in the sense of \cite{NPF05}.
\end{example}
A block substitution naturally determines a collection of patches. Recall that  a patch $P$ is called $S$-legal if there exists a letter $a\in \mathcal{A}$ and $n\in \mathbb{N}$ such that $P\prec S^n(a)$. We denote the set of all $S$-legal patches by $W(S)$, and define the substitution subshift
\begin{equation}
 \Omega(S):= \{ \omega \in \mathcal{A}^{\mathbb{Z}^d}: W(\omega) \subseteq W(S) \}.
\end{equation} 
A patch $P\in W(\omega)\setminus W(S)$ is called a structural defect of $\omega$ (with respect to $S$). Note that $W(S)$ is downwards closed; that is, if $Q\prec P$ and $P\in W(S)$, then $Q\in W(S)$. Moreover, $S$ is closed under the substitution; if $P\in W(S)$, then $S(P)\in W(S)$. 

Following \cite{BBPT24, BBPT26}, we mainly focus on $T$-patches, be they legal or illegal.
Note that 
\[ T \subseteq D_{\mathsf{q}}\big[ T+[0,1)^d \big] \quad \text{and} \quad T^j \subseteq  D_{\mathsf{q}}\big[ T^j+[0,1)^d \big] \]
for every $1\leq j \leq d$.
Following \cite{BBPT26}, we are particularly interested in structural defects arising in fixed points of S. To study such defects, we analyze in Section \ref{Sect:Defect} the ways in which they can occur in fixed points. The connection between fixed points of a substitution and $T$-patches is illustrated in \cite[Proposition 4.3]{BBPT26}.
\begin{proposition} [\cite{BBPT26}] \label{Prop:FixPt}
	Let $S$ be a primitive block substitution of length $\mathsf{q}=(q_1,...,q_d)$, where $q_j\geq 2$ for all $1\leq j \leq d$. Let $Q\in \mathcal{A}^T$ such that $S^p(Q)|_T=Q$ for some $p\in \mathbb{N}$. Then, for every $\rho \in \mathcal{A}^{ \mathbb{Z}^d }$ with $\rho|_T=Q$, the limit
	$$ \rho_Q :=\lim_{n\to \infty} S^{np}(\rho)  $$
	exists in $\mathcal{A}^{ \mathbb{Z}^d }$ and is independent of the choice of $\rho$. Moreover, $\rho_Q|_T= Q$ and $S^p(\rho_Q)=\rho_Q$.
\end{proposition}

\subsection{Number coding of a patch address} \label{Sect:NumCode}
Given \eqref{Eq:SuppIter} and the identity 
$$ D_\mathsf{q}^n\Big[ \operatorname{supp}(P) +[0,1)^d \Big]= D^n_\mathsf{q}\big[ \operatorname{supp}(P) \big]+ D^n_\mathsf{q}\big[ [0,1)^d \big], $$ 
every $\mathsf{m} \in \operatorname{supp}\big( S^n(P) \big)$ is of the form 
$$ \mathsf{m} = \sum_{i=0}^n D_\mathsf{q}^{n-i} (\mathsf{r}^i), \quad \text{where} \quad \mathsf{r}^0\in \supp(P) \;\; \text{and} \;\; \mathsf{r}^i\in K_\mathsf{q}  $$
for all $1\leq i \leq n$. It follows inductively that, for every $1\leq k \leq n$, 
\begin{equation} \label{Eq:SubIter}
\big[ S^k(P) \big](\mathsf{m}^k)= \Big[ S_0\big( \big[S^{k-1}(P)\big](\mathsf{m}^{k-1})\big) \Big](\mathsf{r}^{k}),
\end{equation}
where 
$$ \mathsf{m}^k:= \sum_{i=0}^k D_\mathsf{q}^{k-i} (\mathsf{r}^k) \quad \text{for all} \quad 0\leq k \leq n. $$
For a fixed $\mathsf{q}\in \mathbb{N}^d$, we refer to the tuple $(\mathsf{r}^0;\mathsf{r}^1,...,\mathsf{r}^n)$ as the $n$-th level $\mathsf{q}$-coding, or \Emph{address tuple}, of the vector $\mathsf{m}\in \mathbb{Z}^d$. 
The identity in \eqref{Eq:SubIter} shows that this coding naturally determines an automaton associated with the substitution rule, which in turn uniquely determines the value of a substituted patch at the point $\mathsf{m}$.

We now recall the notions of a finite directed graph and of an automaton.
A \Emph{directed graph} $\mathcal{G}$  consists of a pair $(\mathcal{V},\mathcal{E})$, where $\mathcal{V}$ is the \Emph{vertex set} and $\mathcal{E}\subseteq \mathcal{V}\times \mathcal{V}$ is the \Emph{directed edge set}. 
A tuple $(v_0,v_1,...,v_n)\in \mathcal{V}^{n+1}$ is called a \Emph{walk} of length $n$ starting at $v_0$ and ending at $v_n$, if 
$$ (v_{i-1},v_i)\in \mathcal{E} \quad \text{for every} \quad 1\leq i \leq n. $$
We call the walk \Emph{self-avoiding} if $|\{ v_i \}_{i=1}^n|=n$.
Let $L$ be a finite nonempty set of labels. 
A partially defined function from $X$ to $Y$, denoted $\delta: X 	\rightharpoonup Y$, is a function $\delta:\mathcal{D}\to Y$ defined on some subset $\mathcal{D}\subseteq X$, called its domain.
A partially defined function $\mathcal{\delta}: L\times \mathcal{V} 	\rightharpoonup \mathcal{V}$ satisfying 
$$\big(v, \delta(\ell,v)\big)\in \mathcal{E} \quad \text{whenever} \quad  \delta(\ell,v) \;\; \text{is defined} $$
a \Emph{transition function}. The transition function is called \Emph{complete} if it is defined on all of $L\times \mathcal{V}$, in which case we write $\delta:L\times \mathcal{V}\to \mathcal{V}$. Otherwise, we call it \Emph{incomplete}. 
If $u=\delta(\ell,v)$, we say that $(v,u)$ is \Emph{labelled} by $\ell$. Note that an edge may carry more than one label. We furthermore assume that every edge is labelled; that is, 
\[
(v,u)\in\mathcal E
\quad \text{if and only if} \quad
\exists\,\ell\in L \text{ such that } u=\delta(\ell,v).
\]
We extend $\delta$ recursively to words over $L$. For $\ell'\in L$, $\tilde\ell\in L^n$, and $v\in\mathcal V$, define 
\[ \delta(\ell' \tilde{\ell} ,v):= \delta\big(\ell', \delta(\tilde{\ell},v) \big) \quad \text{whenever both} \quad \delta(\tilde{\ell},v) \; \; \text{and} \;\; \delta\big(\ell', \delta(\tilde{\ell},v) \big) \quad \text{are defined}  . \]
 More generally, the transition function induces transitions labelled by words in $L^*:=\cup_{n=1}L^n$, and therefore induces labels on walks whose lengths coincide with those of the corresponding words.
The triple $(\mathcal{V},\mathcal{E},\delta)$ is called an \Emph{automaton}, complete  or  incomplete according to the transition function.

\begin{defn} \label{Defn:SubAuto}
The substitution rule automaton associated with a $\mathbb{Z}^d$ block substitution over $\mathcal{A}$ is the triple $(\mathcal{A},\mathcal{E},\delta_S)$, where the complete transition function $\delta_S:K_\mathsf{q}\times \mathcal{A}\to \mathcal{A}$ is defined by
 \[ \delta_S(\mathsf{r}, a)= b \quad \text{if and only if} \quad [S_0(a)](\mathsf{r})= b. \]
\end{defn}

For a patch $P$, the value of $S^n(P)$ at a point $\mathsf{m}\in \supp\big( S^n(P) \big)$ can be determined from the substitution rule automaton in the following way. 
If $a_0=P(\mathsf{r}^0)$ and if $a_1,...,a_n$ are defined by
\[ [S^k(P)](\mathsf{m}^k)=a_k \quad \text{for all} \quad k<0 \leq n, \]
then \eqref{Eq:SubIter} can be rewritten as
\[ a_k= \big[ S_0(a_{k-1}) \big](\mathsf{r}^k) \quad \text{for all} \quad 1\leq k\leq n.  \]
This induces a labelled walk in the automaton,
\[ P \overset{\mathsf{r}^0}{\mapsto} a_0 \overset{\mathsf{r}^1 }{\longrightarrow } a_1 \overset{\mathsf{r}^2 }{\longrightarrow } ... \overset{\mathsf{r}^n }{\longrightarrow } a_n .\]
Recall that we assume the dilation vectors $\mathsf{q}$ satisfy $\mathsf{q}_j\geq 2$ for all $1\leq j \leq d$, so that our associated dilations are expanding in all directions. Consequently,
$$ D_\mathsf{q}^n\big[ [0,1)^d \big]= \prod_{j=1}^d \big[0, q_j^n\big)  \to [0,\infty)^d \quad \text{and} \quad  D_\mathsf{q}^n\big[ \{-1,0\}^d+ [0,1)^d \big]= \prod_{j=1}^d \big[ -q_j^n, q_j^n \big)  \to \mathbb{Z}^d $$
monotonically as $n\to \infty$. 
\newline
For $m \in \mathbb{Z}$ and $2\leq q\in \mathbb{N}$, 
define 
\begin{equation} \label{Eq:CodingExpon}
k(q;m):=\min\{ p\in \mathbb{N}: -q^p \leq m< q^p \}.
\end{equation}
Then every $m\in \mathbb{Z}$ admits a unique representation 
\begin{equation} \label{Eq:qCoding}
m= \sum_{i=0}^{k(q;m)} b_i q^{k(q;m)-i},
\end{equation}
where $b_0\in\{-1,0\}$ and $b_i\in \{0,...,q-1\}$ for all $1\leq i \leq k(q;m)$. 
For $0\leq k \leq k(q;m)$, define
\[ m_k:= \sum_{i=0}^k b_i q^{k-i} \quad \text{for all} \quad 0\leq k \leq k(q;m). \] 
We focus on structural defects whose supports are of the form $m\mathsf{e}_j+T^j$ for some $1\leq j \leq d$, for the sets $T^j$ defined in \eqref{Eq:TjDefn}.
Given  $\mathsf{t}\in T^j$ and $m\in \mathbb{N}_0$, with $q_j$-coding $m= \sum_{i=0}^{k(q;m)} b_i q^{k(q;m)-i}$, one has
\begin{equation}
	\mathsf{t}+m \mathsf{e}_j \in D^{k(q_j;m)}\big[ \mathsf{t}+ b_0\mathsf{e}_j+ [0,1)^d \big].
\end{equation}
Moreover, we can uniquely write 
$$m\mathsf{e}_j+ \mathsf{t}=\sum_{i=0}^{k(q_j;m)}D_\mathsf{q}^{k(q;m)-i}(\mathsf{r}^i_{\mathsf{t}}+b_i \mathsf{e}_j),$$ 
where
$\mathsf{r}_{\mathsf{t}}^0=\mathsf{t}$, and 
$$ r_\ell^{(i)}= t_\ell- q_\ell t_\ell \quad \text{for all} \quad 1\leq \ell\leq d, $$
where $\mathsf{r}_\mathsf{t}^i=\big(r_1^{(i)},...,r_d^{(i)}\big)$ and $\mathsf{t}=(t_1,...,t_d)$.
An induction argument using \eqref{Eq:SubIter} shows that there exist sequences of $T^j$-patches $\{ P_k \}_{k=0}^{k(q_j;m)}$ and integers $\{ m_k \}_{k=0}^{k(q_j;m)}$ such that
\begin{equation} \label{Eq:Address1}
 m_k:= \sum_{i=0}^k b_i q_j^{k-i} \quad \text{and} \quad P_k=S^k(P_0)|_{m_k \mathsf{e}_j +T^j} 
\end{equation}
for all $0\leq k\leq k(q_j;m)$, and
\begin{equation} \label{Eq:Address2}
	P_k\prec S(P_{k-1}) \quad \text{for all} \quad 1\leq k \leq k(q_j;m).
\end{equation}
Motivated by the observation in \eqref{Eq:Address1} and \eqref{Eq:Address2}, for every $1\leq j \leq d$ we define the $j$-th axis substitution automaton. Let us first introduce the notion of an induced automaton. Let $(\mathcal{V},\mathcal{E}, \delta)$ be a finite automaton. For a subset $\mathcal{V}'\subseteq \mathcal{V}$, we define the \Emph{induced automaton} $(\mathcal{V'}, \mathcal{E}', \delta')$, where $\mathcal{E}'=\mathcal{E}\cap(\mathcal{V}'\times \mathcal{V}' )$ and 
$$ \delta'(\ell,v)=u \quad \text{if and only if} \quad \delta(\ell,v)=u \;\; \text{and} \;\; (v,u)\in \mathcal{E}'. $$
Note that even if  $\delta$ is a complete transition function, the induced transition function may be incomplete.
\begin{defn} \label{Defn:SubAutoJ}
The \Emph{$j$-th axis substitution automaton} associated with a $\mathbb{Z}^d$ block substitution over $\mathcal{A}$ is the triple $\mathcal{G}_j(S)=(\mathcal{A}^{T^j},\mathcal{E}_j,\delta_j)$, where the complete transition function $\delta_j:\{0,...,q_j-1\}\times \mathcal{A}^{T^j}\to \mathcal{A}^{T^j}$ is defined such that
 \[ \delta_j(r, \check{P})= \check{Q} \quad \text{if and only if} \quad S(\check{P})|_{r\mathsf{e}_j+ T^j}= \check{Q}. \]
The \Emph{$j$-th axis defect automaton}, denoted by $\tilde{\mathcal{G}}_j(S)$, is the induced automaton on $\mathcal{A}^{T^j}\setminus W(S)$.
\end{defn}

It is important to note that the axis substitution automaton always has a complete transition function, we discuss examples in this note where the axis defect automaton has an incomplete transition function.  Let us consider several examples of these automata.
\begin{example} \label{Examp:Squir}
Recall the symbolic squiral substitution given in Figure \ref{Fig:Squiral_Rule} over the binary alphabet $\mathcal{A}=\{0,1\}$. Its dilation map is induced from the vector $(q_1,q_2)=(3,3)$, and the label set we consider is $\{ 0,1,2 \}$. Both its $1$-st and $2$-nd substitution automata are given in Figures \ref{Fig:Squiral-1SubAuto} and \ref{Fig:Squiral-2SubAuto}.


\begin{multicols}{2}
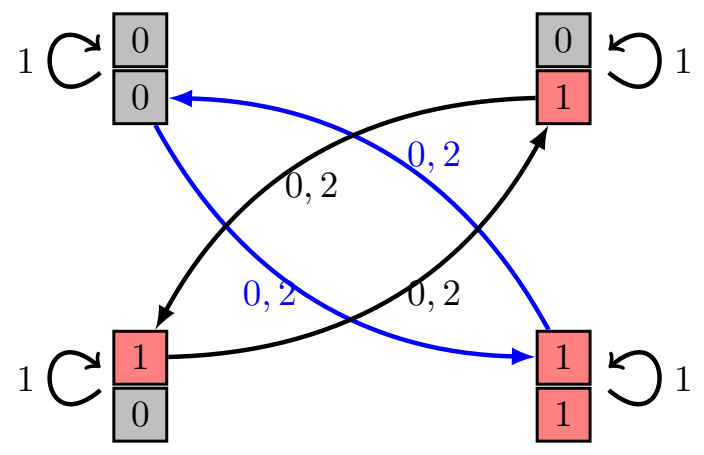
\begin{figure}[H]
	\centering
	\scalebox{0.7}
	{
	\begin{tikzpicture}
				[
				node distance= 0.1mm, 
				0node/.style={rectangle, draw=black, fill=gray!50!white, thick, minimum size=5mm, font ={$0$} },
				1node/.style={rectangle, draw=black, fill=red!50!white, thick, minimum size=5mm, font ={$1$}},
				2node/.style={rectangle, draw=black, fill=green!50!white, thick, minimum size=5mm, font ={$2$}},
				3node/.style={rectangle, draw=black, fill=blue!50!white, thick, minimum size=5mm, font={$3$}},
				]
				
				\begin{scope}[xshift=-2cm, yshift=1.5cm]
				\draw[draw=none,xshift=-0.25cm, yshift=-0.2cm] (0,0) circle  (0.15cm)node (00) {};
				\node[0node] (U00) {};
				\node[0node] (B00) [below= of U00] {};
				\end{scope}
				
				\begin{scope}[xshift=2cm, yshift=1.5cm]
				\draw[draw=none,xshift=0.3cm, yshift=-0.2cm] (0,0) circle  (0.15cm)node (01) {};
				\node[0node] (U01) {};
				\node[1node] (B01) [below= of U01] {};
				\end{scope}
				
				\begin{scope}[xshift=-2cm, yshift=-1.5cm]
				\draw[draw=none,xshift=-0.25cm, yshift=-0.2cm] (0,0) circle  (0.15cm)node (10) {};
				\node[1node] (U10) {};
				\node[0node] (B10) [below= of U10] {};
				\end{scope}
				
				\begin{scope}[xshift=2cm, yshift=-1.5cm]
				\draw[draw=none,xshift=0.3cm, yshift=-0.2cm](0,0) circle  (0.15cm)node (11) {};
				\node[1node] (U11) {};
				\node[1node] (B11) [below= of U11] {};
				\end{scope}
				
				\draw[-latex, ,bend right, very thick, blue] (B00) edge  node[left] {$0,2$}  (U11);	
				\draw[-latex, ,bend right, very thick, blue] (U11) edge  node[right] {$0,2$}  (B00);				
				\draw[-latex, ,bend right, very thick] (U10) edge  node[right] {$0,2$} (B01);	
				\draw[-latex, ,bend right, very thick] (B01) edge  node[below] {$0,2$}  (U10);					
				\draw[->, very thick] (00) to[out=220, in=140, looseness=10, min distance=0.5cm] node[left] {$1$} (00);	
				\draw[->, very thick] (10) to[out=220, in=140, looseness=10, min distance=0.5cm] node[left] {$1$} (10);	
				\draw[->, very thick] (01) to[out=-40, in=40, looseness=10, min distance=0.5cm] node[right] {$1$} (01);	
				\draw[->, very thick] (11) to[out=-40, in=40, looseness=10, min distance=0.5cm] node[right] {$1$} (11);		
				
	\end{tikzpicture}	
	}
	\caption{\centering The $1$-st axis substitution automaton of $S_{\mathtt{Squir}}$.}
	\label{Fig:Squiral-1SubAuto}
\end{figure}

\columnbreak

\begin{figure}[H]
	\centering
	\scalebox{0.7}
	{
	\begin{tikzpicture}
				[
				node distance= 0.1mm, 
				0node/.style={rectangle, draw=black, fill=gray!50!white, thick, minimum size=5mm, font ={$0$} },
				1node/.style={rectangle, draw=black, fill=red!50!white, thick, minimum size=5mm, font ={$1$}},
				2node/.style={rectangle, draw=black, fill=green!50!white, thick, minimum size=5mm, font ={$2$}},
				3node/.style={rectangle, draw=black, fill=blue!50!white, thick, minimum size=5mm, font={$3$}},
				]
				
				\begin{scope}[xshift=-2cm, yshift=1.5cm]
				\draw[draw=none,xshift=-0.35cm, yshift=0cm] (0,0) circle  (0.15cm)node (00) {};
				\node[0node] (L00) {};
				\node[0node] (R00) [right= of L00] {};
				\end{scope}
				
				\begin{scope}[xshift=2cm, yshift=1.5cm]
				\draw[draw=none,xshift=0.85cm, yshift=0cm] (0,0) circle  (0.15cm)node (01) {};
				\node[0node] (L01) {};
				\node[1node] (R01) [right= of L01] {};
				\end{scope}
				
				\begin{scope}[xshift=-2cm, yshift=-1.5cm]
				\draw[draw=none,xshift=-0.35cm, yshift=0cm] (0,0) circle  (0.15cm)node (10) {};
				\node[1node] (L10) {};
				\node[0node] (R10) [right= of L10] {};
				\end{scope}
				
				\begin{scope}[xshift=2cm, yshift=-1.5cm]
				\draw[draw=none,xshift=0.85cm, yshift=0cm](0,0) circle  (0.15cm)node (11) {};
				\node[1node] (L11) {};
				\node[1node] (R11) [right= of L11] {};
				\end{scope}
				
				\draw[-latex, ,bend right, very thick, blue] (R00) edge  node[below left] {$0,2$}  (L11);	
				\draw[-latex, ,bend right, very thick, blue] (L11) edge  node[above right] {$0,2$}  (R00);				
				\draw[-latex, ,bend right, very thick] (R10) edge  node[below right] {$0,2$} (L01);	
				\draw[-latex, ,bend right, very thick] (L01) edge  node[above left] {$0,2$}  (R10);					
				\draw[->, very thick] (00) to[out=220, in=140, looseness=10, min distance=0.5cm] node[left] {$1$} (00);	
				\draw[->, very thick] (10) to[out=220, in=140, looseness=10, min distance=0.5cm] node[left] {$1$} (10);	
				\draw[->, very thick] (01) to[out=-40, in=40, looseness=10, min distance=0.5cm] node[right] {$1$} (01);	
				\draw[->, very thick] (11) to[out=-40, in=40, looseness=10, min distance=0.5cm] node[right] {$1$} (11);		
				
	\end{tikzpicture}
	}
	\caption{\centering The $2$-nd axis substitution automaton of $S_{\mathtt{Squir}}$.}
	\label{Fig:Squiral-2SubAuto}
\end{figure}
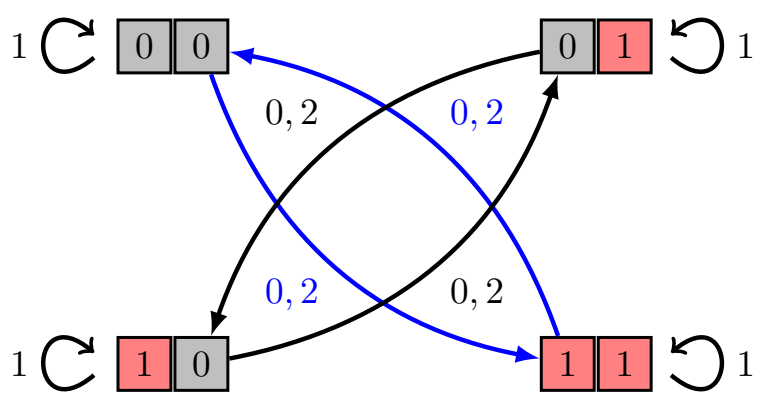

\end{multicols}

Since $W(S)_{T^1}=\mathcal{A}^{T^1}$ and $W(S)_{T^2}=\mathcal{A}^{T^2}$, the $1$-st and $2$-nd axis defect automata are both equal to the empty graph.
\end{example}

\begin{example}
	We turn to observe the axes substitution and defect automata associated with the symbolic table tiling substitution. The $1$-st and $2$-nd axes substitution automata are given in Figures \ref{Fig:Table-1SubAuto} and \ref{Fig:Table-2SubAuto} respectively.

\begin{figure}[h]
	\centering
	\scalebox{0.85}
	{
	\begin{tikzpicture}
		[
				node distance= 0.1mm, 
				0node/.style={rectangle, draw=black, fill=gray!50!white, thick, minimum size=5mm, font ={$0$} },
				1node/.style={rectangle, draw=black, fill=red!50!white, thick, minimum size=5mm, font ={$1$}},
				2node/.style={rectangle, draw=black, fill=green!50!white, thick, minimum size=5mm, font ={$2$}},
				3node/.style={rectangle, draw=black, fill=blue!50!white, thick, minimum size=5mm, font={$3$}},
				]
		
	
	\begin{scope}[xshift=-0.5cm]
		\draw[draw=none,xshift=-0.35cm, yshift=-0.2cm] (0,0) circle  (0.15cm)node (00) {};
		\node[0node] (U00) {};
		\node[0node] (B00) [below= of U00] {};
	\end{scope}
	
	\begin{scope}[yshift=-2cm, xshift=-0.5cm]
		\draw[draw=none,xshift=-0.35cm, yshift=-0.2cm] (0,0) circle  (0.15cm)node (31) {};
		\node[3node] (U31) {};
		\node[1node] (B31) [below= of U31] {};
	\end{scope}
	
	\begin{scope}[yshift=-2cm, xshift=1cm]
		\draw[draw=none,xshift=0.25cm, yshift=-0.2cm] (0,0) circle  (0.15cm)node (22) {};
		\node[2node] (U22) {};
		\node[2node] (B22) [below= of U22] {};
	\end{scope}
	
	\begin{scope}[xshift=0.75cm]
		\draw[draw=none,xshift=0.25cm, yshift=-0.2cm] (0,0) circle  (0.15cm)node (13) {};
		\node[1node] (U13) {};
		\node[3node] (B13) [below= of U13] {};
	\end{scope}
	
	
	\begin{scope}[xshift=3cm]
		\draw[draw=none,xshift=-0.35cm, yshift=-0.2cm] (0,0) circle  (0.15cm)node (01) {};
		\node[0node] (U01) {};
		\node[1node] (B01) [below= of U01] {};
	\end{scope}
	
	\begin{scope}[xshift=4.5cm]
		\draw[draw=none,xshift=-0.35cm, yshift=-0.2cm] (0,0) circle  (0.15cm)node (30) {};
		\node[3node] (U30) {};
		\node[0node] (B30) [below= of U30] {};
	\end{scope}
	
	\begin{scope}[xshift=6cm]
		\draw[draw=none,xshift=-0.35cm, yshift=-0.2cm] (0,0) circle  (0.15cm)node (20) {};
		\node[2node] (U20) {};
		\node[0node] (B20) [below= of U20] {};
	\end{scope}
	
	\begin{scope}[yshift=-2cm, xshift=3cm]
		\draw[draw=none,xshift=-0.35cm, yshift=-0.2cm] (0,0) circle  (0.15cm)node (02) {};
		\node[0node] (U02) {};
		\node[2node] (B02) [below= of U02] {};
	\end{scope}
	
	\begin{scope}[yshift=-2cm, xshift=4.5cm]
		\draw[draw=none,xshift=-0.35cm, yshift=-0.2cm] (0,0) circle  (0.15cm)node (32) {};
		\node[3node] (U32) {};
		\node[2node] (B32) [below= of U32] {};
	\end{scope}
	
	\begin{scope}[yshift=-2cm, xshift=6cm]
		\draw[draw=none,xshift=-0.35cm, yshift=-0.2cm] (0,0) circle  (0.15cm)node (21) {};
		\node[2node] (U21) {};
		\node[1node] (B21) [below= of U21] {};
	\end{scope}

	\draw[-latex, ,bend right, very thick] (B00) edge  node[left] {$0$} (U31);	
	\draw[-latex, ,bend right, very thick] (U31) edge node[right] {$0$} (B00);		
	
	\draw[-latex, ,bend right, very thick] (U22) edge  node[below] {$1$} (U31);	
	\draw[-latex, ,bend right, very thick] (B31) edge node[above] {$1$} (B22);

	\draw[-latex, ,bend right, very thick] (U30) edge  node[below] {$0$} (U01);	
				\draw[-latex, ,bend right, very thick] (B01) edge node[above] {$0$}   (B30);	
				\draw[-latex, ,bend right, very thick] (U20) edge  node[below] {$1$}  (U30);	
				\draw[-latex, ,bend right, very thick] (B30) edge node[above] {$1$} (B20);	
				\draw[-latex, ,bend right, very thick] (U32) edge  node[below] {$0$} (U02);	
				\draw[-latex, ,bend right, very thick] (B02) edge node[above] {$0$} (B32);	
				\draw[-latex, ,bend right, very thick] (U21) edge node[below] {$1$}   (U32);	
				\draw[-latex, ,bend right, very thick] (B32) edge node[above] {$1$} (B21);	
				
	\draw[-latex, ,bend right, very thick] (B01) edge  node[left] {$1$} (U02);	
	\draw[-latex, ,bend right, very thick] (U02) edge node[right] {$1$} (B01);	
	
	\draw[-latex, ,bend right, very thick] (B20) edge  node[left] {$0$} (U21);	
	\draw[-latex, ,bend right, very thick] (U21) edge node[right] {$0$} (B20);	
	
	\draw[->, very thick] (00) to[out=220, in=140, looseness=10, min distance=0.5cm] node[left] {$1$} (00);	
	\draw[->, very thick] (13) to[out=-40, in=40, looseness=10, min distance=0.5cm] node[right] {$0,1$} (13);	
	\draw[->, very thick] (22) to[out=-40, in=40, looseness=10, min distance=0.5cm] node[right] {$0$} (22);

	
	\begin{scope}[xshift=9.5cm, yshift=-1cm]
	
		\begin{scope}[xshift=-1.5cm, yshift=1cm]
				\draw[draw=none,xshift=-0.35cm, yshift=-0.2cm] (0,0) circle  (0.15cm)node (10) {};
				\node[1node] (U10) {};
				\node[0node] (B10) [below= of U10] {};
				\end{scope}	
				
				\begin{scope}[xshift=0cm, yshift=1cm]
				\node[1node] (U11) {};
				\node[1node] (B11) [below= of U11] {};
				\end{scope}				
				
				\begin{scope}[xshift=1.5cm, yshift=1cm]
				\draw[draw=none,xshift=0.35cm, yshift=-0.2cm] (0,0) circle  (0.15cm)node (12) {};
				\node[1node] (U12) {};
				\node[2node] (B12) [below= of U12] {};
				\end{scope}	
				
				\begin{scope}[xshift=-1.5cm, yshift=-1cm]
				\draw[draw=none,xshift=-0.35cm, yshift=-0.2cm] (0,0) circle  (0.15cm)node (03) {};
				\node[0node] (U03) {};
				\node[3node] (B03) [below= of U03] {};
				\end{scope}	
				
				\begin{scope}[xshift=0cm, yshift=-1cm]
				\node[3node] (U33) {};
				\node[3node] (B33) [below= of U33] {};
				\end{scope}				
				
				\begin{scope}[xshift=1.5cm, yshift=-1cm]
				\draw[draw=none,xshift=0.35cm, yshift=-0.2cm] (0,0) circle  (0.15cm)node (23) {};
				\node[2node] (U23) {};
				\node[3node] (B23) [below= of U23] {};
				\end{scope}	
	
	\end{scope}
	
	\draw[-latex, ,bend right, very thick] (U11) edge  node[below] {$0$} (U10);	
				\draw[-latex, ,bend right, very thick] (B10) edge node[above] {$0$}   (B11);	
				\draw[-latex, ,bend right, very thick] (U12) edge  node[below] {$1$}  (U11);	
				\draw[-latex, ,bend right, very thick] (B11) edge node[above] {$1$} (B12);	
				\draw[-latex, ,bend right, very thick] (U33) edge  node[below] {$0$} (U03);	
				\draw[-latex, ,bend right, very thick] (B03) edge node[above] {$0$} (B33);	
				\draw[-latex, ,bend right, very thick] (U23) edge node[below] {$1$}   (U33);	
				\draw[-latex, ,bend right, very thick] (B33) edge node[above] {$1$} (B23);	
				\draw[->, very thick] (10) to[out=220, in=140, looseness=10, min distance=0.5cm] node[left] {$1$} (10);	
				\draw[->, very thick] (03) to[out=220, in=140, looseness=10, min distance=0.5cm] node[left] {$1$} (03);	
				\draw[->, very thick] (12) to[out=-40, in=40, looseness=10, min distance=0.5cm] node[right] {$0$} (12);	
				\draw[->, very thick] (23) to[out=-40, in=40, looseness=10, min distance=0.5cm] node[right] {$0$} (23);	
		
	\end{tikzpicture}
	}
	\caption{\centering The $1$-st axis substitution automaton of $S_{\mathtt{Table}}$.}
	\label{Fig:Table-1SubAuto}
\end{figure}
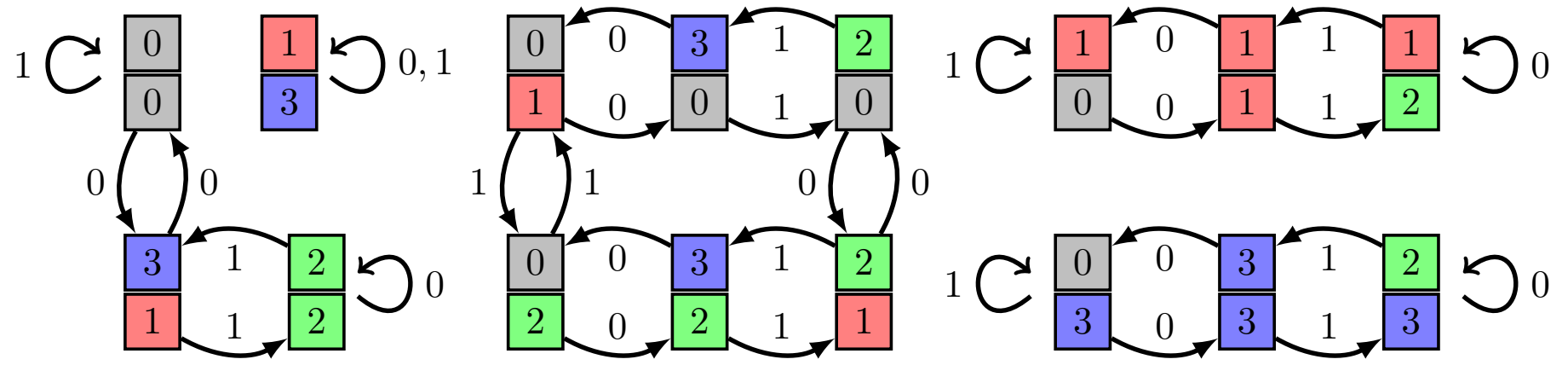

\begin{figure}[h]
	\centering
	\scalebox{0.85}
	{
	\begin{tikzpicture}
		[
				node distance= 0.1mm, 
				0node/.style={rectangle, draw=black, fill=gray!50!white, thick, minimum size=5mm, font ={$0$} },
				1node/.style={rectangle, draw=black, fill=red!50!white, thick, minimum size=5mm, font ={$1$}},
				2node/.style={rectangle, draw=black, fill=green!50!white, thick, minimum size=5mm, font ={$2$}},
				3node/.style={rectangle, draw=black, fill=blue!50!white, thick, minimum size=5mm, font={$3$}},
				scale=0.9
				]
				
	
	\begin{scope}[yshift=1cm]
			\draw[draw=none,xshift=0.25cm, yshift=-0.2cm] (0,0) circle  (0.15cm)node (02) {};
				\node[0node] (L02) {};
				\node[2node] (R02) [right= of L02] {};
	\end{scope}
	
	\begin{scope}[yshift=-1cm]
			\draw[draw=none,xshift=-0.25cm, yshift=0cm] (0,0) circle  (0.15cm)node (11) {};
				\node[1node] (L11) {};
				\node[1node] (R11) [right= of L11] {};
	\end{scope}
	
	\begin{scope}[yshift=1cm, xshift =2.25cm]
			\draw[draw=none,xshift=-0.25cm, yshift=0cm] (0,0) circle  (0.15cm)node (33) {};
				\node[3node] (L33) {};
				\node[3node] (R33) [right= of L33] {};
	\end{scope}
	
	\begin{scope}[yshift=-1cm, xshift=2cm]
			\draw[draw=none,xshift=-0.25cm, yshift=0cm] (0,0) circle  (0.15cm)node (20) {};
				\node[2node] (L20) {};
				\node[0node] (R20) [right= of L20] {};
	\end{scope}
	
	\begin{scope}[yshift=1cm, xshift=4cm]
			\draw[draw=none,xshift=-0.25cm, yshift=0cm] (0,0) circle  (0.15cm)node (10) {};
				\node[1node] (L10) {};
				\node[0node] (R10) [right= of L10] {};
	\end{scope}
	
	\begin{scope}[yshift=-1cm, xshift=4cm]
			\draw[draw=none,xshift=-0.25cm, yshift=0cm] (0,0) circle  (0.15cm)node (13) {};
				\node[1node] (L13) {};
				\node[3node] (R13) [right= of L13] {};
	\end{scope}
	
	\begin{scope}[yshift=1cm, xshift=6cm]
			\draw[draw=none,xshift=-0.25cm, yshift=0cm] (0,0) circle  (0.15cm)node (21) {};
				\node[2node] (L21) {};
				\node[1node] (R21) [right= of L21] {};
	\end{scope}
	
	\begin{scope}[yshift=-1cm, xshift=6cm]
			\draw[draw=none,xshift=-0.25cm, yshift=0cm] (0,0) circle  (0.15cm)node (23) {};
				\node[2node] (L23) {};
				\node[3node] (R23) [right= of L23] {};
	\end{scope}
	
	\begin{scope}[yshift=1cm, xshift=8cm]
			\draw[draw=none,xshift=-0.25cm, yshift=0cm] (0,0) circle  (0.15cm)node (31) {};
				\node[3node] (L31) {};
				\node[1node] (R31) [right= of L31] {};
	\end{scope}
	
	\begin{scope}[yshift=-1cm, xshift=8cm]
			\draw[draw=none,xshift=-0.25cm, yshift=0cm] (0,0) circle  (0.15cm)node (30) {};
				\node[3node] (L30) {};
				\node[0node] (R30) [right= of L30] {};
	\end{scope}
	
	\draw[-latex, ,bend right, very thick] (L21) edge  node[above] {$1$} (R10);	
				\draw[-latex, ,bend right, very thick] (R10) edge node[below] {$1$} (L21);	
	
	\draw[-latex, ,bend right, very thick] (L31) edge  node[above] {$0$} (R21);	
				\draw[-latex, ,bend right, very thick] (R21) edge node[below] {$0$} (L31);	
				
	\draw[-latex, ,bend right, very thick] (L23) edge  node[above] {$1$} (R13);	
				\draw[-latex, ,bend right, very thick] (R13) edge node[below] {$1$} (L23);	
				
	\draw[-latex, ,bend right, very thick] (L30) edge  node[above] {$0$} (R23);	
				\draw[-latex, ,bend right, very thick] (R23) edge node[below] {$0$} (L30);	
				
	\draw[-latex, ,bend right, very thick] (R13) edge  node[left] {$0$} (R10);	
				\draw[-latex, ,bend right, very thick] (L10) edge node[right] {$0$} (L13);	
				
	\draw[-latex, ,bend right, very thick] (R30) edge  node[left] {$1$} (R31);	
				\draw[-latex, ,bend right, very thick] (L31) edge node[right] {$1$} (L30);

	\draw[-latex, ,bend right, very thick] (L20) edge  node[above] {$1$} (R11);	
				\draw[-latex, ,bend right, very thick] (R11) edge node[below] {$1$} (L20);	
				
	\draw[-latex, ,bend right, very thick] (R20) edge  node[left] {$0$} (R33);	
				\draw[-latex, ,bend right, very thick] (L33) edge node[right] {$0$} (L20);
				
	\draw[->, very thick] (11) to[out=220, in=140, looseness=10, min distance=0.5cm] node[left] {$0$} (11);
	\draw[->, very thick] (33) to[out=220, in=140, looseness=10, min distance=0.5cm] node[left] {$0$} (33);	
	\draw[->, very thick] (02) to[out=310, in=230, looseness=10, min distance=0.5cm] node[below] {$0,1$} (02);

	\begin{scope}[xshift=12.5cm, yshift=0cm]
		\begin{scope}[xshift=-2cm, yshift=1cm]
				\draw[draw=none,xshift=-0.25cm, yshift=0cm] (0,0) circle  (0.15cm)node (01) {};
				\node[0node] (L01) {};
				\node[1node] (R01) [right= of L01] {};
				\end{scope}	
				
				\begin{scope}[xshift=0cm, yshift=1cm]
				\node[0node] (L00) {};
				\node[0node] (R00) [right= of L00] {};
				\end{scope}				
				
				\begin{scope}[xshift=2cm, yshift=1cm]
				\draw[draw=none,xshift=0.85cm, yshift=0cm] (0,0) circle  (0.15cm)node (03) {};
				\node[0node] (L03) {};
				\node[3node] (R03) [right= of L03] {};
				\end{scope}	
				
				\begin{scope}[xshift=-2cm, yshift=-1cm]
				\draw[draw=none,xshift=-0.25cm, yshift=0cm] (0,0) circle  (0.15cm)node (12) {};
				\node[1node] (L12) {};
				\node[2node] (R12) [right= of L12] {};
				\end{scope}	
				
				\begin{scope}[xshift=0cm, yshift=-1cm]
				\node[2node] (L22) {};
				\node[2node] (R22) [right= of L22] {};
				\end{scope}				
				
				\begin{scope}[xshift=2cm, yshift=-1cm]
				\draw[draw=none,xshift=0.85cm, yshift=0cm] (0,0) circle  (0.15cm)node (32) {};
				\node[3node] (L32) {};
				\node[2node] (R32) [right= of L32] {};
				\end{scope}
				
				\draw[-latex, ,bend right, very thick] (L00) edge  node[above] {$1$} (R01);	
				\draw[-latex, ,bend right, very thick] (R01) edge node[below] {$1$} (L00);	
				\draw[-latex, ,bend right, very thick] (L03) edge  node[above] {$0$} (R00);	
				\draw[-latex, ,bend right, very thick] (R00) edge node[below] {$0$} (L03);	
				\draw[-latex, ,bend right, very thick] (R12) edge  node[below] {$1$} (L22);	
				\draw[-latex, ,bend right, very thick] (L22) edge node[above] {$1$}  (R12);	
				\draw[-latex, ,bend right, very thick] (R22) edge  node[below] {$0$}  (L32);	
				\draw[-latex, ,bend right, very thick] (L32) edge node[above] {$0$} (R22);	
				
				\draw[->, very thick] (01) to[out=220, in=140, looseness=10, min distance=0.5cm] node[left] {$0$} (01);	
				\draw[->, very thick] (12) to[out=220, in=140, looseness=10, min distance=0.5cm] node[left] {$0$} (12);	
				\draw[->, very thick] (03) to[out=-40, in=40, looseness=10, min distance=0.5cm] node[right] {$1$} (03);	
				\draw[->, very thick] (32) to[out=-40, in=40, looseness=10, min distance=0.5cm] node[right] {$1$} (32);	
	\end{scope}			
				
	\end{tikzpicture}
	}
	\caption{\centering The $2$-nd axis substitution automaton of $S_{\mathtt{Table}}$.}
	\label{Fig:Table-2SubAuto}
\end{figure}
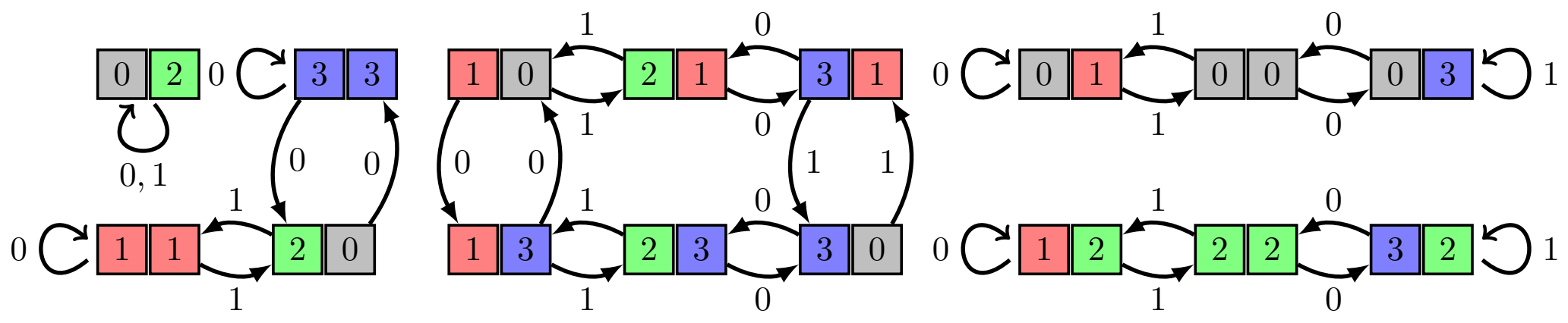

We can restrict these automata to the vertex set of illegal patches and obtain Figures \ref{Fig:Table-1DefectAuto} and \ref{Fig:Table-2DefectAuto}.


\begin{multicols}{2}
	
\begin{figure}[H]
\centering
	\scalebox{0.85}
	{
	\begin{tikzpicture}
				[
				node distance= 0.1mm, 
				0node/.style={rectangle, draw=black, fill=gray!50!white, thick, minimum size=5mm, font ={$0$} },
				1node/.style={rectangle, draw=black, fill=red!50!white, thick, minimum size=5mm, font ={$1$}},
				2node/.style={rectangle, draw=black, fill=green!50!white, thick, minimum size=5mm, font ={$2$}},
				3node/.style={rectangle, draw=black, fill=blue!50!white, thick, minimum size=5mm, font={$3$}},
				]
				
				\begin{scope}[xshift=-2cm, yshift=1cm]
				\draw[draw=none,xshift=-0.35cm, yshift=-0.2cm] (0,0) circle  (0.15cm)node (10) {};
				\node[1node] (U10) {};
				\node[0node] (B10) [below= of U10] {};
				\end{scope}	
				
				\begin{scope}[xshift=0cm, yshift=1cm]
				\node[1node] (U11) {};
				\node[1node] (B11) [below= of U11] {};
				\end{scope}				
				
				\begin{scope}[xshift=2cm, yshift=1cm]
				\draw[draw=none,xshift=0.35cm, yshift=-0.2cm] (0,0) circle  (0.15cm)node (12) {};
				\node[1node] (U12) {};
				\node[2node] (B12) [below= of U12] {};
				\end{scope}	
				
				\begin{scope}[xshift=-2cm, yshift=-1cm]
				\draw[draw=none,xshift=-0.35cm, yshift=-0.2cm] (0,0) circle  (0.15cm)node (03) {};
				\node[0node] (U03) {};
				\node[3node] (B03) [below= of U03] {};
				\end{scope}	
				
				\begin{scope}[xshift=0cm, yshift=-1cm]
				\node[3node] (U33) {};
				\node[3node] (B33) [below= of U33] {};
				\end{scope}				
				
				\begin{scope}[xshift=2cm, yshift=-1cm]
				\draw[draw=none,xshift=0.35cm, yshift=-0.2cm] (0,0) circle  (0.15cm)node (23) {};
				\node[2node] (U23) {};
				\node[3node] (B23) [below= of U23] {};
				\end{scope}
				
				\draw[-latex, ,bend right, very thick] (U11) edge  node[below] {$0$} (U10);	
				\draw[-latex, ,bend right, very thick] (B10) edge node[above] {$0$}   (B11);	
				\draw[-latex, ,bend right, very thick] (U12) edge  node[below] {$1$}  (U11);	
				\draw[-latex, ,bend right, very thick] (B11) edge node[above] {$1$} (B12);	
				\draw[-latex, ,bend right, very thick] (U33) edge  node[below] {$0$} (U03);	
				\draw[-latex, ,bend right, very thick] (B03) edge node[above] {$0$} (B33);	
				\draw[-latex, ,bend right, very thick] (U23) edge node[below] {$1$}   (U33);	
				\draw[-latex, ,bend right, very thick] (B33) edge node[above] {$1$} (B23);	
				\draw[->, very thick] (10) to[out=220, in=140, looseness=10, min distance=0.5cm] node[left] {$1$} (10);	
				\draw[->, very thick] (03) to[out=220, in=140, looseness=10, min distance=0.5cm] node[left] {$1$} (03);	
				\draw[->, very thick] (12) to[out=-40, in=40, looseness=10, min distance=0.5cm] node[right] {$0$} (12);	
				\draw[->, very thick] (23) to[out=-40, in=40, looseness=10, min distance=0.5cm] node[right] {$0$} (23);	
				
	\end{tikzpicture}
	}
	\caption{\centering The $1$-st axis defect automaton of $S_{\mathtt{Table}}$.}
	\label{Fig:Table-1DefectAuto}
\end{figure}
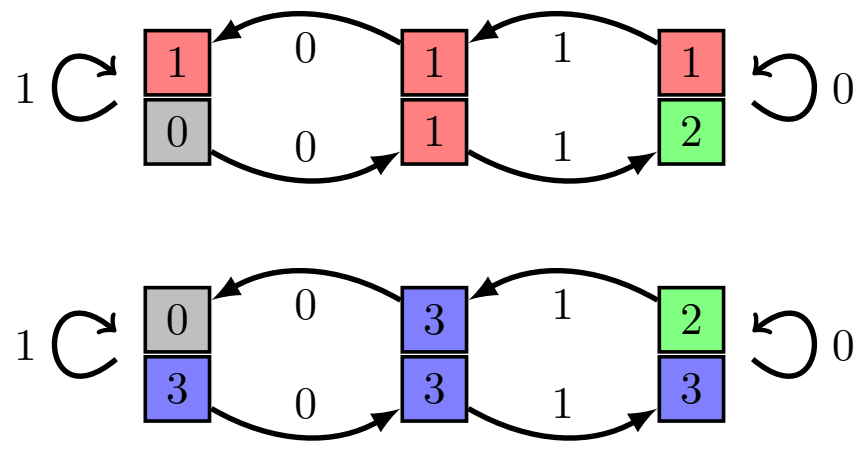	
	
\columnbreak

	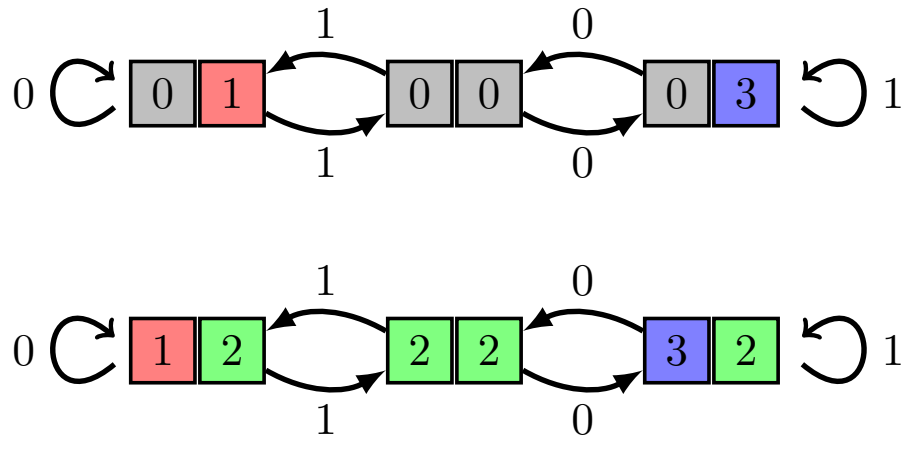
\begin{figure}[H]
\centering
	\scalebox{0.85}
	{
	\begin{tikzpicture}
				[
				node distance= 0.1mm, 
				0node/.style={rectangle, draw=black, fill=gray!50!white, thick, minimum size=5mm, font ={$0$} },
				1node/.style={rectangle, draw=black, fill=red!50!white, thick, minimum size=5mm, font ={$1$}},
				2node/.style={rectangle, draw=black, fill=green!50!white, thick, minimum size=5mm, font ={$2$}},
				3node/.style={rectangle, draw=black, fill=blue!50!white, thick, minimum size=5mm, font={$3$}},
				]
				
				\begin{scope}[xshift=-2cm, yshift=1cm]
				\draw[draw=none,xshift=-0.25cm, yshift=0cm] (0,0) circle  (0.15cm)node (01) {};
				\node[0node] (L01) {};
				\node[1node] (R01) [right= of L01] {};
				\end{scope}	
				
				\begin{scope}[xshift=0cm, yshift=1cm]
				\node[0node] (L00) {};
				\node[0node] (R00) [right= of L00] {};
				\end{scope}				
				
				\begin{scope}[xshift=2cm, yshift=1cm]
				\draw[draw=none,xshift=0.85cm, yshift=0cm] (0,0) circle  (0.15cm)node (03) {};
				\node[0node] (L03) {};
				\node[3node] (R03) [right= of L03] {};
				\end{scope}	
				
				\begin{scope}[xshift=-2cm, yshift=-1cm]
				\draw[draw=none,xshift=-0.25cm, yshift=0cm] (0,0) circle  (0.15cm)node (12) {};
				\node[1node] (L12) {};
				\node[2node] (R12) [right= of L12] {};
				\end{scope}	
				
				\begin{scope}[xshift=0cm, yshift=-1cm]
				\node[2node] (L22) {};
				\node[2node] (R22) [right= of L22] {};
				\end{scope}				
				
				\begin{scope}[xshift=2cm, yshift=-1cm]
				\draw[draw=none,xshift=0.85cm, yshift=0cm] (0,0) circle  (0.15cm)node (32) {};
				\node[3node] (L32) {};
				\node[2node] (R32) [right= of L32] {};
				\end{scope}
				
				\draw[-latex, ,bend right, very thick] (L00) edge  node[above] {$1$} (R01);	
				\draw[-latex, ,bend right, very thick] (R01) edge node[below] {$1$} (L00);	
				\draw[-latex, ,bend right, very thick] (L03) edge  node[above] {$0$} (R00);	
				\draw[-latex, ,bend right, very thick] (R00) edge node[below] {$0$} (L03);	
				\draw[-latex, ,bend right, very thick] (R12) edge  node[below] {$1$} (L22);	
				\draw[-latex, ,bend right, very thick] (L22) edge node[above] {$1$}  (R12);	
				\draw[-latex, ,bend right, very thick] (R22) edge  node[below] {$0$}  (L32);	
				\draw[-latex, ,bend right, very thick] (L32) edge node[above] {$0$} (R22);	
				
				\draw[->, very thick] (01) to[out=220, in=140, looseness=10, min distance=0.5cm] node[left] {$0$} (01);	
				\draw[->, very thick] (12) to[out=220, in=140, looseness=10, min distance=0.5cm] node[left] {$0$} (12);	
				\draw[->, very thick] (03) to[out=-40, in=40, looseness=10, min distance=0.5cm] node[right] {$1$} (03);	
				\draw[->, very thick] (32) to[out=-40, in=40, looseness=10, min distance=0.5cm] node[right] {$1$} (32);	
				
	\end{tikzpicture}
	}
	\caption{\centering The $2$-nd axis defect automaton of $S_{\mathtt{Table}}$.}
	\label{Fig:Table-2DefectAuto}
\end{figure}

\end{multicols}

\end{example}

\section{Formation of structural defects via defect automata} \label{Sect:Defect}
	In this section, we describe mechanisms governing the emergence of structural defects in fixed points. 
	Let $S$ be a block substitution with a fixed point $\rho\in \mathcal{A}^{\mathbb{Z}^d}$ satisfying $S^p(\rho)=\rho$ and fix $1\leq j \leq d$. Following Subsection \ref{Sect:NumCode}, every  number $m\in \mathbb{Z}\setminus \{-1,0\}$ admits a unique address $(b_0;b_1,...,b_{k(q_j;m)})$ as in \eqref{Eq:qCoding}. Moreover, there exist $\ell \in \mathbb{N}$ and $\hat{r}\in \mathbb{N}_0$ such that $k(q_j;m)+\hat{r}=\ell p$ and $0\leq \hat{r}<q_j$.
	If $\hat{r}=0$, then $k(q_j;m)$ is a multiple of $p$, and therefore $S^{k(q_j;m)}(\rho)=\rho$. 
	In this case, the induced labelled walk in the $j$-th axis substitution automaton,
	\begin{equation} \label{Eq:FixPtAddress}
	\rho \overset{b_0}{\mapsto} \rho|_{b_0\mathsf{e}_j+T^j} \overset{b_1}{\longrightarrow} S(\rho)|_{(b_0q_j+b_1)\mathsf{e}_j+T^j} \overset{b_2}{\longrightarrow} ... \overset{b_{k(q_j;m)}}{\longrightarrow} S^{k(q_j;m)}(\rho)|_{m\mathsf{e}_j+T^j},
	\end{equation}
	together with \eqref{Eq:Address1}, shows that $\rho|_{m\mathsf{e}_j+T^j}= S^{k(q_j;m)}(\rho)|_{m\mathsf{e}_j+T^j}$. Suppose now that $\hat{r}\neq 0$. Since $\hat{r}+k(q_j;m)$ is a multiple of $p$, we have $S^{\hat{r}+k(q_j;m)}(\rho)= \rho$.  We therefore consider the modified address $$(b_0;\underbrace{b_0(1-q_j),...,b_0(1-q_j)}_{\hat{r}\; \text{times}},b_1,...,b_{k(q_j;m)}),$$
	which introduces a delay of $\hat{r}$ substitution steps. 
	Using \eqref{Eq:Address1}, we conclude that the  patch $\rho|_{m\mathsf{e}_j+T^j}= S^{\hat{r}+k(q_j;m)}(\rho)|_{m\mathsf{e}_j+T^j}$ is determined by the corresponding labelled walk in the j-th axis substitution automaton,
	\begin{equation} \label{Eq:GenFixPtAddress}
	\begin{aligned}
	\rho \overset{b_0}{\mapsto} \rho|_{b_0\mathsf{e}_j+T^j} &\overset{b_0(1-q_j)}{\longrightarrow} S(\rho)|_{b_0\mathsf{e}_j+T^j} \overset{b_0(1-q_j)}{\longrightarrow} ... \\
	&\overset{b_0(1-q_j)}{\longrightarrow} S^{\hat{r}}(\rho)|_{b_0\mathsf{e}_j+T^j} \overset{b_1}{\longrightarrow} ... \overset{b_{k(q_j;m)}}{\longrightarrow} S^{k(q_j;m)}(\rho)|_{m\mathsf{e}_j+T^j}.
	\end{aligned}
	\end{equation}

The following lemma summarizes the implications of these addresses for fixed points of block substitutions.

\begin{lemma} \label{Lem:DefctAuto}
	Let $S$ be a block substitution of length $\mathsf{q}=(q_1,...,q_d)$ and $\rho\in \mathcal{A}^{\mathbb{Z}^d}$ be a fixed point of $S$, with $S^p(\rho)=\rho$ for some $p\in \mathbb{N}$. Fix some $1\leq j \leq d$ and let $m\in \mathbb{Z}\setminus \{-1,0\}$ with the $q_j$-number coding in \eqref{Eq:qCoding} given by
	\[ m= \sum^{k(q_j;m)}_{i=0} b_i q_j^{k(q_j;m)-i}. \]
	Then, the following statements hold.
	\begin{enumerate}[(a)]
		\item There exist $\ell\in \mathbb{N}$ and $\hat{r}\in \mathbb{N}_0$, such that $k(q_j;m)+\hat{r}= \ell p$ and $0\leq\hat{r}<p$. 
		The patch $\rho|_{m\mathsf{e}_j+T^j}$ is the endpoint of the walk labelled $(\underbrace{b_0(1-q_j),...,b_0(1-q_j)}_{\hat{r}\; \text{times}},b_1,...,b_{k(q_j;m)})$ starting at $\rho|_{b_0\mathsf{e}_j+T^j}$, as described in \eqref{Eq:GenFixPtAddress}.
		\item The patch $\rho|_{m\mathsf{e}_j+T^j}$ is a structural defect if and only if the labelled walk described in (a) also defines a walk in the $j$-th axis defect automaton with the same label. 
	\end{enumerate}
	
\end{lemma}	
	
Let us consider the following example to illustrate this mechanism.
\begin{example} \label{Examp:DefectForm}
Consider the symbolic table tiling substitution.
Following Proposition \ref{Prop:FixPt}, we know that a fixed point $\rho$ of $S_{\mathtt{Table}}$ is determined uniquely by the patch $\rho|_{ \{-1,0 \}^2 }$. Consider the illegal $T$-patch  $P_0:\{-1,0\}^2\to \{0,...,3\}$ in Figure \ref{Fig:P0Patch}.

\begin{multicols}{2}

\begin{figure}[H]
	\centering
	\begin{tikzpicture}
		[
				node distance= 0.1mm, 
				0node/.style={rectangle, draw=black, fill=gray!50!white, thick, minimum size=5mm, font ={$0$} },
				1node/.style={rectangle, draw=black, fill=red!50!white, thick, minimum size=5mm, font ={$1$}},
				2node/.style={rectangle, draw=black, fill=green!50!white, thick, minimum size=5mm, font ={$2$}},
				3node/.style={rectangle, draw=black, fill=blue!50!white, thick, minimum size=5mm, font={$3$}},
				]
				
				\node[0node] (LU) {};
				\node[1node] (RU) [right= of LU] {};
				\node[3node] (LB) [below= of LU] {};
				\node[2node] [right= of LB] {};
				
	\end{tikzpicture}
	\caption{The patch $P_0$.}
	\label{Fig:P0Patch}
\end{figure}


\columnbreak


\begin{figure}[H]
	\centering
	\begin{tikzpicture}
		[
				node distance= 0.1mm, 
				0node/.style={rectangle, draw=black, fill=gray!50!white, thick, minimum size=5mm, font ={$0$} },
				1node/.style={rectangle, draw=black, fill=red!50!white, thick, minimum size=5mm, font ={$1$}},
				2node/.style={rectangle, draw=black, fill=green!50!white, thick, minimum size=5mm, font ={$2$}},
				3node/.style={rectangle, draw=black, fill=blue!50!white, thick, minimum size=5mm, font={$3$}},
				]
				
				\node[1node] (LU) {};
				\node[1node] (RU) [right= of LU] {};
				\node[1node] (LB) [below= of LU] {};
				\node[1node] [right= of LB] {};
				
	\end{tikzpicture}
	\caption{The patch $P_1$.}
	\label{Fig:P1Patch}
\end{figure}

\end{multicols}
The patch $P_0$ induces a unique fixed point $S_{\mathtt{Table}}(\rho_{P_0})=\rho_{P_0}$. We can therefore use the derivation in \eqref{Eq:FixPtAddress} to determine a $T^j$-structural defect in $\rho_{P_0}$. Let us take $m\in \mathbb{Z}$ such that $k(2;m)=\min\{ p\in \mathbb{N}:-2^p\leq m < 2^p  \}=3$. Take $m'=6$ and $m''=-6$, which satisfy
\[ 6 = 0\cdot 2^3 + 1\cdot 2^2 + 1\cdot 2^1+ 0 \cdot 2^0 \quad \text{and} \quad -6= -1\cdot 2^3 +0\cdot 2^2 +1\cdot 2^1 +0 \cdot 2^0, \]
which implies that $(0;1,1,0)$ and $(-1;0,1,0)$ are the $2$-codings of $m'=6$ and $m''-6$ respectively. Since $P_0|_{T^1}, P_0|_{-\mathsf{e}_1+T^1}$ are not $S_{\mathtt{Table}}$-legal, we can use the $1$-st axis
defect automaton of $S_{\mathtt{Table}}$ to determine ${\rho_{P_0}|}_{m'' \mathsf{e}_1 +T^1}$. 
The appropriate derivation following \eqref{Eq:FixPtAddress} is given in Figure \ref{Fig:m''Deriv}, using the $1$-st axis defect automaton in Figure \ref{Fig:Table-1DefectAuto}.

\begin{multicols}{2}

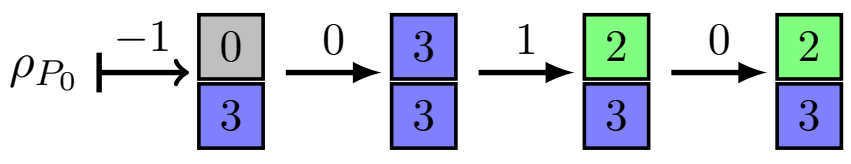
\begin{figure}[H]
	\centering
	\scalebox{0.85}
	{
	\begin{tikzpicture}
		[
				node distance= 0.1mm, 
				0node/.style={rectangle, draw=black, fill=gray!50!white, thick, minimum size=5mm, font ={$0$} },
				1node/.style={rectangle, draw=black, fill=red!50!white, thick, minimum size=5mm, font ={$1$}},
				2node/.style={rectangle, draw=black, fill=green!50!white, thick, minimum size=5mm, font ={$2$}},
				3node/.style={rectangle, draw=black, fill=blue!50!white, thick, minimum size=5mm, font={$3$}},
				]
				
			\begin{scope}[xshift=-1.75cm]
			\draw[draw=none,xshift=0.3cm, yshift=-0.2cm] (0,0) circle  (0.15cm)node (R-1) {\large $\rho_{P_0}$};
			\end{scope}				
				
			\begin{scope}[xshift=0cm]
			\draw[draw=none,xshift=-0.2cm, yshift=-0.2cm] (0,0) circle  (0.15cm)node (L0) {};
			\draw[draw=none,xshift=0.3cm, yshift=-0.2cm] (0,0) circle  (0.15cm)node (R0) {};
				\node[0node] (U0) {};
				\node[3node] (B0) [below= of U0] {};
			\end{scope}
			
			\begin{scope}[xshift=1.5cm]
			\draw[draw=none,xshift=-0.2cm, yshift=-0.2cm] (0,0) circle  (0.15cm)node (L1) {};
			\draw[draw=none,xshift=0.3cm, yshift=-0.2cm] (0,0) circle  (0.15cm)node (R1) {};
				\node[3node] (U1) {};
				\node[3node] (B1) [below= of U1] {};
			\end{scope}	
			
			\begin{scope}[xshift=3cm]
			\draw[draw=none,xshift=-0.2cm, yshift=-0.2cm] (0,0) circle  (0.15cm)node (L2) {};
			\draw[draw=none,xshift=0.3cm, yshift=-0.2cm] (0,0) circle  (0.15cm)node (R2) {};
				\node[2node] (U2) {};
				\node[3node] (B2) [below= of U2] {};
			\end{scope}	
			
			\begin{scope}[xshift=4.5cm]
			\draw[draw=none,xshift=-0.2cm, yshift=-0.2cm] (0,0) circle  (0.15cm)node (L3) {};
			\draw[draw=none,xshift=0.3cm, yshift=-0.2cm] (0,0) circle  (0.15cm)node (R3) {};
				\node[2node] (U3) {};
				\node[3node] (B3) [below= of U3] {};
			\end{scope}	
				
			\draw[|->, very thick] (R-1) edge node[above] {$-1$} (L0);		
			\draw[-latex, very thick] (R0) edge node[above] {$0$} (L1);				
			\draw[-latex, very thick] (R1) edge node[above] {$1$} (L2);
			\draw[-latex, very thick] (R2) edge node[above] {$0$} (L3);		
	\end{tikzpicture}
	}
	\caption{\centering The derivation for $\rho_{P_0}|_{m''\mathsf{e}_1+T^1}$.}
	\label{Fig:m''Deriv}
\end{figure}

\columnbreak

\begin{figure}[H]
	\centering
	\scalebox{0.85}
	{
	\begin{tikzpicture}
		[
				node distance= 0.1mm, 
				0node/.style={rectangle, draw=black, fill=gray!50!white, thick, minimum size=5mm, font ={$0$} },
				1node/.style={rectangle, draw=black, fill=red!50!white, thick, minimum size=5mm, font ={$1$}},
				2node/.style={rectangle, draw=black, fill=green!50!white, thick, minimum size=5mm, font ={$2$}},
				3node/.style={rectangle, draw=black, fill=blue!50!white, thick, minimum size=5mm, font={$3$}},
				scale=0.9
				]
				
			\begin{scope}[xshift=-1.75cm]
			\draw[draw=none,xshift=0.3cm, yshift=-0.2cm] (0,0) circle  (0.15cm)node (R-1) {\large $\rho_{P_1}$};
			\end{scope}				
				
			\begin{scope}[xshift=0cm]
			\draw[draw=none,xshift=-0.2cm, yshift=-0.2cm] (0,0) circle  (0.15cm)node (L0) {};
			\draw[draw=none,xshift=0.3cm, yshift=-0.2cm] (0,0) circle  (0.15cm)node (R0) {};
				\node[1node] (U0) {};
				\node[1node] (B0) [below= of U0] {};
			\end{scope}
			
			\begin{scope}[xshift=1.5cm]
			\draw[draw=none,xshift=-0.2cm, yshift=-0.2cm] (0,0) circle  (0.15cm)node (L1) {};
			\draw[draw=none,xshift=0.3cm, yshift=-0.2cm] (0,0) circle  (0.15cm)node (R1) {};
				\node[1node] (U1) {};
				\node[0node] (B1) [below= of U1] {};
			\end{scope}	
			
			\begin{scope}[xshift=3cm]
			\draw[draw=none,xshift=-0.2cm, yshift=-0.2cm] (0,0) circle  (0.15cm)node (L2) {};
			\draw[draw=none,xshift=0.3cm, yshift=-0.2cm] (0,0) circle  (0.15cm)node (R2) {};
				\node[1node] (U2) {};
				\node[0node] (B2) [below= of U2] {};
			\end{scope}	
			
			\begin{scope}[xshift=4.5cm]
			\draw[draw=none,xshift=-0.2cm, yshift=-0.2cm] (0,0) circle  (0.15cm)node (L3) {};
			\draw[draw=none,xshift=0.3cm, yshift=-0.2cm] (0,0) circle  (0.15cm)node (R3) {};
				\node[1node] (U3) {};
				\node[0node] (B3) [below= of U3] {};
			\end{scope}	
			
			\begin{scope}[xshift=6cm]
			\draw[draw=none,xshift=-0.2cm, yshift=-0.2cm] (0,0) circle  (0.15cm)node (L4) {};
			\draw[draw=none,xshift=0.3cm, yshift=-0.2cm] (0,0) circle  (0.15cm)node (R4) {};
				\node[1node] (U4) {};
				\node[1node] (B4) [below= of U4] {};
			\end{scope}
				
			\draw[|->, very thick] (R-1) edge node[above] {$0$} (L0);		
			\draw[-latex, very thick] (R0) edge node[above] {$0$} (L1);				
			\draw[-latex, very thick] (R1) edge node[above] {$1$} (L2);
			\draw[-latex, very thick] (R2) edge node[above] {$1$} (L3);
			\draw[-latex, very thick] (R3) edge node[above] {$0$} (L4);		
	\end{tikzpicture}
	}
	\caption{\centering The derivation for $\rho_{P_1}|_{m'\mathsf{e}_1+T^1}$.}
	\label{Fig:m'DerivGen}
\end{figure}
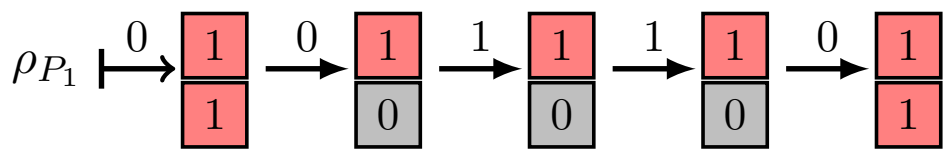

\end{multicols}

\begin{figure}[b]
	\centering
	\begin{tikzpicture}
				[
				node distance= 0.1mm, 
				0node/.style={rectangle, draw=black, fill=gray!50!white, thick, minimum size=5mm, font ={$0$} },
				1node/.style={rectangle, draw=black, fill=red!50!white, thick, minimum size=5mm, font ={$1$}},
				2node/.style={rectangle, draw=black, fill=green!50!white, thick, minimum size=5mm, font ={$2$}},
				3node/.style={rectangle, draw=black, fill=blue!50!white, thick, minimum size=5mm, font={$3$}},
				]
				\node[1node] (0U) {};
				\node[2node] (0B) [below= of 0U]  {};
				\node[0node] (-1U) [left= of 0U]  {};
				\node[3node] (-1B) [below= of -1U]  {};
				
				\node[1node] (1U) [right= of 0U]  {};
				\node[1node] (1B) [below= of 1U]  {};
				\node[1node] (2U) [right= of 1U]  {};
				\node[0node] (2B) [below= of 2U]  {};
				\node[1node] (3U) [right= of 2U]  {};
				\node[2node] (3B) [below= of 3U]  {};
				
				\node[1node] (4U) [right= of 3U]  {};
				\node[1node] (4B) [below= of 4U]  {};
				\node[1node] (5U) [right= of 4U]  {};
				\node[0node] (5B) [below= of 5U]  {};
				\node[1node] (6U) [right= of 5U]  {};
				\node[2node] (6B) [below= of 6U]  {};
				
				\node[1node] (7U) [right= of 6U]  {};
				\node[1node] (7B) [below= of 7U]  {};
				\node[1node] (8U) [right= of 7U]  {};
				\node[0node] (8B) [below= of 8U]  {};
				\node[1node] (9U) [right= of 8U]  {};
				\node[2node] (9B) [below= of 9U]  {};
				
				\node[1node] (10U) [right= of 9U]  {};
				\node[1node] (10B) [below= of 10U]  {};
				\node[1node] (11U) [right= of 10U]  {};
				\node[0node] (11B) [below= of 11U]  {};
				\node[1node] (12U) [right= of 11U]  {};
				\node[2node] (12B) [below= of 12U]  {};	
				
				\node[1node] (13U) [right= of 12U]  {};
				\node[1node] (13B) [below= of 13U]  {};
				
				\node[3node] (-2U) [left= of -1U]  {};
				\node[3node] (-2B) [below= of -2U]  {};
				\node[2node] (-3U) [left= of -2U]  {};
				\node[3node] (-3B) [below= of -3U]  {};
				\node[0node] (-4U) [left= of -3U]  {};
				\node[3node] (-4B) [below= of -4U]  {};	
				
				\node[3node] (-5U) [left= of -4U]  {};
				\node[3node] (-5B) [below= of -5U]  {};
				\node[2node] (-6U) [left= of -5U]  {};
				\node[3node] (-6B) [below= of -6U]  {};
				\node[0node] (-7U) [left= of -6U]  {};
				\node[3node] (-7B) [below= of -7U]  {};	
				
				\node[3node] (-8U) [left= of -7U]  {};
				\node[3node] (-8B) [below= of -8U]  {};
				\node[2node] (-9U) [left= of -8U]  {};
				\node[3node] (-9B) [below= of -9U]  {};
				\node[0node] (-10U) [left= of -9U]  {};
				\node[3node] (-10B) [below= of -10U]  {};	
				
				\node[3node] (-11U) [left= of -10U]  {};
				\node[3node] (-11B) [below= of -11U]  {};
				\node[2node] (-12U) [left= of -11U]  {};
				\node[3node] (-12B) [below= of -12U]  {};
				\node[0node] (-13U) [left= of -12U]  {};
				\node[3node] (-13B) [below= of -13U]  {};
				
				\node[3node] (-14U) [left= of -13U]  {};
				\node[3node] (-14B) [below= of -14U]  {};

				\draw[orange, fill=yellow, fill opacity=0.35,  xshift=-0.75cm, yshift=-0.25cm] (0,-0.55) rectangle (0.5,0.5);
				
				\draw[orange, fill=yellow, fill opacity=0.35,  xshift=-1.28cm, yshift=-0.25cm] (0,-0.55) rectangle (0.5,0.5);
				
				\draw[orange, fill=yellow, fill opacity=0.35,  xshift=-1.85cm, yshift=-0.25cm] (0,-0.55) rectangle (0.5,0.5);
				0,-0.5) rectangle (0.5,0.5);
				
				\draw[orange, fill=yellow, fill opacity=0.35,  xshift=-3.45cm, yshift=-0.25cm] (0,-0.55) rectangle (0.5,0.5);
				
				
	\end{tikzpicture}
	\caption{\centering A portion of the pure-line defect supported on $\Lambda_1$ in the fixed point $\rho_{P_0}$, with the patches in Figure \ref{Fig:m''Deriv}.}
	\label{Fig:PureLineDefectSimpl}
\end{figure}
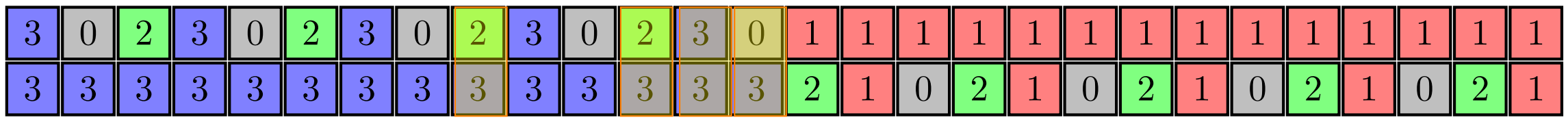

To illustrate how to study fixed points satisfying $S^2_{\mathtt{Table}}(\rho)=\rho$, we take the fixed point coming from the constant illegal $T$-patch $P_1(t)\equiv 1$, shown in Figure \ref{Fig:P1Patch}. By Proposition \ref{Prop:FixPt}, the induced fixed point $\rho_{P_1}$ satisfies $S^2_{\mathtt{Table}}(\rho_{P_1})= \rho_{P_1}$. To determine ${\rho_{P_1}}|_{m'\mathsf{e}_1+T^1}$ using the $1$-st axis defect automaton, as described in \eqref{Eq:GenFixPtAddress},we first note that the modified address corresponding to $m'=6$ and $P_1$ is $(0;0,1,1,0)$ with a delay of $\hat{r}=1$ step.  The corresponding derivation is presented in Figure \ref{Fig:m'DerivGen}, while Figure \ref{Fig:PureLineDefectConst1} shows a portion of the pure-line defect occurring in $\rho_{P_1}$.

\end{example}

\begin{proof}[Proof of Lemma \ref{Lem:DefctAuto}]
$~$\\
 (a) 
Let $\ell= \big\lceil \frac{k(q_j;m)}{p} \big\rceil$ and $\hat{r}= \ell p-k(q_j;m)$. Note that by our assumptions, $k:=k(q_j;m)\geq 1$.
By \eqref{Eq:SuppIter}, $m\mathsf{e}_j+T^j \subseteq \supp\big( S^{k}(\rho|_{b_0\mathsf{e}_j+ T^j}) \big)$, where $b_0=0$ if $m\geq 0$ and $b_0=-1$ otherwise. Since $S^p(\rho)=\rho$, we have $S^{k+\hat{r}}(\rho)|_{m\mathsf{e}_j+T^j}=\rho|_{m\mathsf{e}_j+T^j}$. Set $P_0:= S^{\hat{r}}(\rho)|_{b_0\mathsf{e}_j+T^j}$. By \eqref{Eq:Address1} and \eqref{Eq:Address2}, there exist patches $P_1,...,P_{k}$ such that $P_{k-1} \overset{b_k}{\longrightarrow}P_k$ and $P_{k}=S^{k}(P_0)|_{m \mathsf{e}_j+T^j}$. Consequently,
\[ \rho|_{-\mathsf{e_j}+T^j} \overset{b_0(1-q_j)}{\longrightarrow} S(\rho)|_{-\mathsf{e_j}+T^j} \overset{b_0(1-q_j)}{\longrightarrow} ... \overset{b_0(1-q_j)}{\longrightarrow}  S^{\hat{r}}(\rho)|_{-\mathsf{e_j}+T^j} \overset{b_1}{\longrightarrow} P_1 \overset{b_2}{\longrightarrow} ... \overset{b_{k}}{\longrightarrow} P_{k} = \rho|_{m\mathsf{e}_j+T^j},  \]
is a walk in the j-th axis substitution automaton labelled by $$(\underbrace{b_0(1-q_j),...,b_0(1-q_j)}_{\hat{r}\; \text{times}},b_1,...,b_{k}),$$ as claimed.
\newline
(b) Recall that $W(S)$ is a set of patches that is downward closed with respect to the subpatch relation and closed under the action of the substitution.
Since $P:=\rho|_{m\mathsf{e}_j+T^j}$ is a structural defect, it follows that if $P\prec S(Q)$, for some $T^j$-patch $Q$, then $Q\notin W(S)$. 

In particular, every vertex appearing in the above labelled walk is itself a structural defect. Therefore, the labelled walk is entirely contained in the j-th axis defect automaton, and thus defines a labelled walk in that automaton.
\end{proof}

\subsection{Finitely many bounded structural defects}

	We turn to see what conditions generate only finitely many bounded structural defects, similar to the symbolic chair tiling substitution.	
	The following simple useful result of Lemma \ref{Lem:DefctAuto} satisfies such a condition, for the applicability of \cite[Theorem 2.2]{BBPT26}.

\begin{Cor} \label{Cor:BddDefect}
	Let $S$ be a primitive block substitution of length $\mathsf{q}=(q_1,...,q_d)$, $\rho \in \mathcal{A}^{\mathbb{Z}^d}$ be a fixed point of $S$. If there exists an $N\in \mathbb{N}_0$ such that for all $1\leq j \leq d$, every walk in $\tilde{\mathcal{G}}_j(S)$ has a length less than $N$, then every patch $P\in W(\rho)$ with $\supp(P)\subseteq \mathsf{m}+B_r$, must satisfy that $\{-1,0\}^d \subseteq \supp(P)$ or $P\in W(S)$. 
\end{Cor}

\begin{remark}
An explicit analysis of the defect automata for the symbolic chair tiling substitution shows that the assumptions of Corollary~\ref{Cor:BddDefect} are satisfied with $N=1$.
The corollary can therefore be used to verify the hypotheses of \cite[Theorem~2.2]{BBPT26} for substitution systems beyond the symbolic chair tiling. 
For instance, the symbolic squiral substitution from Example \ref{Examp:Squir} satisfies the assumptions of Corollary~\ref{Cor:BddDefect}, since $\{0,1\}^{T^1}\subseteq W(S_{\mathtt{Squir}})$ and $\{0,1\}^{T^2}\subseteq W(S_{\mathtt{Squir}})$. Hence, the conclusions of \cite[Corollary~2.3(b)]{BBPT26} apply to it as well.
\end{remark}

\begin{proof}
	This is proved using the same argument as \cite[Lemma 5.9]{BBPT26}. For all $1\leq j \leq d$, note that 
	\begin{equation}
	S^n(\rho)|_{T^j} \in W(S) \quad \text{and} \quad S^n(\rho)|_{-\mathsf{e}_j+T^j}\in W(S) 
	\end{equation}
	for $n>N$. 
	Let $P\in W(\rho)$ and assume that $\{-1,0\}^d \not\subseteq \supp(P)$. Hence, there must exist some $1\leq j_0\leq d$ such that 
	$$\mathsf{m}+B_r\subseteq \{ (x_1,...,x_d)\in \mathbb{R}^d: x_{j_0}\geq0\} \quad \text{or} \quad \mathsf{m}+B_r\subseteq \{ (x_1,...,x_d)\in \mathbb{R}^d: x_{j_0}<0\}. $$	
	This implies the existence of an $n_0\in \mathbb{N}$ such that $\mathsf{m}+B_r \subseteq D_{\mathsf{q}}^{n_0}\big[ T^{j_0}+[0,1)^d \big]$. Since $W(S)$ is closed-downwards, with respect to the subpatch relation, \eqref{Eq:SuppIter} implies that $P\in W(S)$.
\end{proof}

\subsection{pure-line defects}

We turn to discuss conditions to give more general conditions for pure-line defects in fixed points. 
Recall that a fixed point $\rho \in \mathcal{A}^{\mathbb{Z}^d}$ of a substitution has a pure-line defect if there exists some $1\leq j \leq d$ such that every $T$-patch occurring in
\[ \Lambda_j:= \{ (m_1,...,m_d)\in \mathbb{Z}^d: m_j\in \mathbb{Z}, \; m_i \in \{-1,0\}\; \text{for}\; i\neq j \} = \mathbb{Z}\mathsf{e}_j+T \]
is $S$-illegal.
Recall from Proposition \ref{Prop:FixPt}, that a fixed point $\rho$ of a substitution is uniquely determined by the value of $P=\rho|_{ \{-1,0\}^d }$.

\begin{lemma} \label{Lem:PL-Defect}
	Let $S$ be a primitive block substitution of length $\mathsf{q}=(q_1,...,q_d)$, and let $\rho \in \mathcal{A}^{\mathbb{Z}^d}$ be a fixed point of $S$ satisfying $S^p(\rho)=\rho$. Fix $1\leq j \leq d$, and set $P_+:= \rho|_{T^j}$ and $P_-:= \rho|_{-\mathsf{e}_j+T^j}$. 
	\begin{enumerate}[(a)]
	\item The fixed point $\rho$ has infinitely many $T^j$-structural defects if and only if there exist distinct labels
$\ell_1,\ell_2\in\{0,\ldots,q_j-1\}^p$
and a patch $P\in\{P_-,P_+\}$ such that both $\ell_1$ and $\ell_2$ label walks in the $j$-th defect automaton starting at $P$.
	\item The fixed point $\rho$ has a pure-line defect along the $j$-th axis if and only if, for every label
$\ell\in\{0,\ldots,q_j-1\}^{N}$, where $N=|\mathcal{A}|^{2^{d-1}}+1$, the walks starting at  $P_+$ and $P_-$ labelled by $\ell$ are  walks in the $j$-th defect automaton.
	\end{enumerate}		
\end{lemma}

Lemma~\ref{Lem:PL-Defect}~(a) provides a necessary and sufficient condition for a fixed point to have infinitely many $T^j$-structural defects, while Lemma~\ref{Lem:PL-Defect}~(b) provides a necessary and sufficient condition for the fixed point to admit a pure-line defect along the $j$-th axis.

\begin{proof}$~$\\
(a)
Since $S^p(\rho)=\rho$, the walks labelled 
$$(\underbrace{0,...,0}_{p \; \text{times}  }) \quad \text{and} \quad (\underbrace{q_j-1,...,q_j-1}_{p \; \text{times}  })$$ starting in $P^+$ and $P^-$ respectively, are both closed walks in $\mathcal{G}_j(S)$. Since each vertex in the walks must be illegal, they are also walks in $\tilde{\mathcal{G}}_j(S)$, the $j$-th axis defect automaton.  

Assume that there exist two distinct labelled walks in the $j$-th defect automaton starting at either $P_+$ or $P_-$, with labels $\ell^{(1)}, \ell^{(2)}\in \{0,...,q_j-1 \}^p$. Without loss of generality, assume that both walks start at $P_+$. Write $\ell^{(i)}=(b_1^{(i)},...,b_p^{(i)})$ for $i\in \{ 1,2 \}$. Since $\ell^{(1)}\neq \ell^{(2)}$, we know that
$$ \sum_{s=1}^p b_s^{(1)}q_j^{p-s}\neq \sum_{s=1}^p b_s^{(2)}q_j^{p-s}.  $$
Without loss of generality, assume that $\sum_{s=0}^p b_s^{(1)}q_j^{p-s}< \sum_{s=0}^p b_s^{(2)}q_j^{p-s}$. For every $k\in \mathbb{N}$, set
$$ m_k= \sum_{i=1}^k\sum_{s=1}^{p} b_s^{(2)} q_j^{s+ip}, $$
and note that $m_k< m_{k+1}$ for all $k\in \mathbb{N}$. By Lemma \ref{Lem:DefctAuto},  $\rho|_{m_k\mathsf{e}_j+T^j}$ is a structural defect for all $k$. Hence, the fixed point $\rho$ has infinitely many $T^j$-structural defects.

Conversely, assume that the fixed point has infinitely many $T^j$-structural defects. There exists an $m\in \mathbb{Z}\setminus \{-1,0\}$ such that $\rho|_{m\mathsf{e}_j+T^j}$ is a structural defect and $m\geq q_j^p$ or $m<-q_j^p$. Assume without of loss of generality that $m<0$. Consider the $q_j$-number coding of $m$ from \eqref{Eq:qCoding},
\[ m= \sum^{k}_{i=0} b_i q_j^{k(q_j;m)-i}, \quad \text{with} \quad k=k(q_j;m). \]
Then $b_0=-1$, and there exists $1\leq i_0= \min\{ i: b_i\neq q_j-1 \}$.
Moreover, $i_0\leq p$. Indeed,
$$ m= -q_j^k +\sum_{i=1}^{i_0-1}(q_j-1)q_j^{k-i} + \sum_{i_0}^{k} b_i q_j^{k-i}\geq -q_j^{k-i_0},  $$
whereas we assumed that $m<-q^{k-1}$.
Since $i_0\leq p$ and $k\geq p$, the labels $(q_j-1,...,q_j-1),(b_0,...,b_p)\in \{ 0,...,q_j-1 \}^p$ are both distinct elements corresponding to walks in $\tilde{\mathcal{G}_j}(S)$ starting at $P_-$.

(b) By Lemma \ref{Lem:DefctAuto}, for every $m\in \mathbb{Z}$ the patch $\rho|_{m\mathsf{e}_j+T^j}$ is a structural defect if and only if for every $n\in \mathbb{N}$ and every label $\ell \in \{0,...,q_j-1\}^n$, $\ell$ is a label for walks in the $j$-th defect automaton starting at $P_-$ and $P_+$.
Hence, it immediately follows that if $\rho$ has a pure-line defect along the $j$-th axis it implies that for every $P\in \{ P_-,P_+ \}$ there exists a walk in $ \tilde{\mathcal{G}}_j(S)$ starting at $P$ and labelled by $\ell=\ell_1...\ell_N\in \{ 0,...,q_j-1 \}^N$. 

Assume now that for every $P\in \{ P_-,P_+ \}$ there exists a walk in $\tilde{\mathcal{G}}_j(S)$ starting at $P$ and labelled by $\ell=\ell_1...\ell_N\in \{ 0,...,q_j-1 \}^N$. 
Assume towards contraposition that $\rho$ does not have a pure-line defect along the $j$-th axis. Hence, there exists by Lemma \ref{Lem:DefctAuto} an $n\in \mathbb{N}$ and a  label $\ell_1...\ell_n\in \{ 0,...,q_j-1\}^n$ such that $\delta(\ell_1...\ell_n,P)\in W(S)$. 
The label induces the following walk in the $j$-th substitution automaton
\[ P \overset{\ell_1}{\longrightarrow} \underbrace{\delta_j(\ell_1,P)}_{P_1} \overset{\ell_2}{\longrightarrow} ... \overset{\ell_n}{\longrightarrow} \underbrace{\delta_j(\ell_1...\ell_n,P)}_{P_n}. \]
The walk induces a self-avoiding walk in the $j$-th substitution automaton $(P,P'_1,...,P'_{n'})$, where $n'\leq n$. By Definition \ref{Defn:SubAutoJ}, there exist labels $\ell_1',...,\ell'_{n'}\in \{ 0,...,q_j-1 \}$ such that 
\[ P'_k = \delta_j(\ell'_1...\ell'_k, P) \quad \text{for all} \quad 1\leq k \leq n'. \]
The pigeonhole principle implies a bound on the length of the induced self-avoiding walk of  $n' \leq |\mathcal{A}^{T^j}|+1 = N $.  Consider the extended label
$$ \ell':=\ell'_1...\ell'_{n'}\hspace*{-0.55cm} \underbrace{0...0}_{(N-n')\; \text{-times}} \hspace*{-0.5cm} \in \{0 ,...,q_j-1\}^N.  $$
By the assumption, there is a walk in $\tilde{\mathcal{G}}_j(S)$ starting at $P$ and labelled by $\ell' \in \{ 0,...,q_j-1 \}^N$. 
Since every walk labelled by $\ell'$ in $\tilde{\mathcal{G}}_j(S)$ starting at $P$ is also a walk labelled by $\ell'$ in $\mathcal{G}_j(S)$, Lemma \ref{Lem:DefctAuto} implies that $\delta_j(\ell', P)$ is a contained in a walk of $\tilde{\mathcal{G}}_j(S)$. This contradicts our contraposition assumption that $\delta_j(\ell', P)= \delta_j(\ell_1...\ell_n, P)\in W(S)$.
\end{proof}

While Lemma~\ref{Lem:PL-Defect}~(b) provides a necessary and sufficient combinatorial condition for the hypothesis of \cite[Proposition~2.8]{BBPT26}, the following corollary gives a simpler sufficient condition for the existence of pure-line defects.

\begin{Cor} \label{Cor:PL-Defect}
Let $1\leq j\leq d$ and let $S$ be a block substitution of length $\mathsf{q}=(q_1,...,q_d)$ such that the $j$-th axis defect automaton $\tilde{G}_j(S)$ is complete. 
Then, any fixed point $\rho\in \mathcal{A}^{\mathbb{Z}^d}$, satisfying $\rho|_{T^j}\notin W(S)$ and $\rho|_{-\mathsf{e}_j+T^j}\notin W(S)$, has a pure-line defect along the $j$-th axis. 
\end{Cor}

\begin{remark}
	Observing Figures \ref{Fig:Table-1DefectAuto} and \ref{Fig:Table-2DefectAuto}, we see that the symbolic table tiling substitution satisfies the assumptions of Corollary~\ref{Cor:PL-Defect} for any $j\in \{1,2\}$. The assumptions of Corollary~\ref{Cor:PL-Defect} are also much easier to check than those of Lemma~\ref{Lem:PL-Defect}~(b).
\end{remark}

\begin{proof}
	Let $\rho\in \mathcal{A}^{\mathbb{Z}^d}$ be a fixed point such that $\rho|_{T^j}\notin W(S)$ and $\rho|_{-\mathsf{e}_j+T^j}\notin W(S)$.
	As in the proof of Lemma~\ref{Lem:PL-Defect}, it follows from Lemma~\ref{Lem:DefctAuto}(b) that $\rho$ has a pure-line defect along the $j$-th axis provided that, for every label $\tilde{\ell}=\ell_1...\ell_n\in \{ 0, ..., q_j-1\}^n $ and every patch $P\in \mathcal{A}^{T^j}\setminus W(S)$, there exists a walk in the $j$-th defect automaton starting at $P$ and labelled by $\tilde{\ell}$. 
	Since $\tilde{\mathcal{G}}_j(S)$ is complete, such a walk exists for every label $\tilde{\ell}$ and every starting vertex $P$. The conclusion therefore follows from Lemma~\ref{Lem:DefctAuto}~(b).
\end{proof}
We turn to discuss a class of substitutions which satisfies the assumptions of Corollary~\ref{Cor:PL-Defect} and naturally gives rise to pure-line defects.

\begin{defn} \label{Defn:Bij}
	Let $S$ be a block substitution of length $\mathsf{q}=(q_1,...,q_d)$.  We say that $S$ is \Emph{bijective along the $j$-th axis}, if for every 
	$$ \mathsf{t}\in K_{\mathsf{q}} \cap  \big\{ (x_1,...,x_d)\in \mathbb{Z}^d :x_j\in \{0 ,q_{j}-1\} \big\},$$ the component map $S_0^{\mathsf{t}}:\mathcal{A} \to \mathcal{A}$ defined by $[S_0^{\mathsf{t}}](a):= [S_0(a)](\mathsf{t})$ is bijective.
\end{defn}

Note that if $S$ is bijective in the sense of \cite{NPF05}, then it is bijective along the $j$-th axis for every $1\leq j \leq d$. Moreover, the two notions coincide when $q_j=2$. 
Bijective block substitutions have been studied extensively because of the special spectral properties of their associated Koopman operators; see, for example, \cite{NPF05,BG14,Bart18}.

The following proposition shows that this class of substitutions is a natural source of pure-line defects.

\begin{proposition} \label{Prop:PureLineDefect}
	Let $S$ be a block substitution of length $\mathsf{q}=(q_1,...,q_d)$. Assume that there exists a $1\leq j \leq d$ such that $S$ is bijective along the $j$-th axis. 
	Then, any fixed point $\rho\in \mathcal{A}^{\mathbb{Z}}$ with $\rho|_{T^j},\rho|_{-\mathsf{e}_j+T^j}\notin W(S)$ has a pure-line defect along the $j$-th axis.
\end{proposition}

\begin{remarks}$~$\
\begin{itemize}
\item Recall that the symbolic table tiling substitution is bijective, and hence is bijective along both the $1$-st and the $2$-nd axes. Proposition \ref{Prop:PureLineDefect} therefore provides a systematic way to determine which fixed points admit pure-line defects, in contrast to the ad hoc arguments used in \cite{BBPT26}. 
\item Bijectivity alone does not guarantee the existence of a fixed point with a pure-line defect along the $j$-th axis. Indeed, if the patches at $T^j$ and $-\mathsf{e}_j+T^j$ are legal, then no pure-line defect along the $j$-th axis can occur. The symbolic squiral substitution provides an example of this phenomenon.
\end{itemize}

\end{remarks}

Using Corollary~\ref{Cor:PL-Defect}, Proposition~\ref{Prop:PureLineDefect} reduces to showing that $\tilde{\mathcal{G}}_j(S)$ is complete whenever $S$ is bijective along the $j$-th axis. We postpone the proof of Proposition~\ref{Prop:PureLineDefect} until Section~\ref{Sect:Tensor}, after introducing the notion of tensor products of automata.

\section{Graph theoretical simplifications}
\label{Sect:Simpl}

\subsection{Graph tensor product structure} \label{Sect:Tensor}
	The automata introduced in Definition~\ref{Defn:SubAutoJ} also admit a tensor-product structure, analogous to that described in \cite[Section~6]{BBPT26}. This decomposition expresses the associated graphs in terms of smaller and more manageable components, thereby reducing the computational complexity of their analysis. Consequently, many of the results and proof techniques from \cite[Section~6]{BBPT26} carry over, with suitable modifications, to the present setting.
	Recall that given directed graphs $\{ G_s=(\mathcal{V}_s,\mathcal{E}_s)  \}_{s=1}^N$, their \Emph{graph tensor product} $\otimes_{s=1}^N G_S$ is the directed graph with vertex set $\prod_{s=1}^N \mathcal{V}_s$ and 
\[ \big( (a_1,...,a_N), (b_1,...b_N) \big)\in \mathcal{E}_{ \otimes_{s=1}^N \mathcal{E}_s } \quad \text{if and only if} \quad (a_s,b_s)\in \mathcal{E}_s \quad \text{for all} \; 1\leq s\leq N.  \]	
If moreover, there exist complete transitions functions $\{\delta_{s}:L\times \mathcal{V} \to \mathcal{V}_s\}_{s=1}^N$, then we have the \Emph{tensor product transition function}
 $\delta_{\otimes_{s=1}^N}:L\times\prod_{s=1}^N \mathcal{V}_s \to  \prod_{s=1}^N \mathcal{V}_s$ satisfying
 \[ \delta_{\otimes_{s=1}^N}\big( (a_1,...,a_N), \ell \big)= (b_1,...,b_N) \quad \text{if and only if} \quad \delta_s(\ell,a_s) = b_s \quad \text{for all} \; 1\leq s\leq N.  \]
The induced determinstic automaton $\big(\otimes_{s=1}^N G_S,  \mathcal{E}_{ \otimes_{s=1}^N \mathcal{E}_s }, \delta_{\otimes_{s=1}^N} \big)$ is called their \Emph{tensor product automaton}.
Fix some $1\leq j \leq d$, and recall that the set $T^j$ is given by
\[ T^j= \{ -1,0 \}^d \cap \{(x_1,...,x_d): x_j=0\}. \]
\begin{defn} \label{Defn:TensorAuto}
For every $\mathsf{t}\in T^j$, define the $(\mathsf{t},j)$-substitution automaton $\mathcal{G}_{\mathsf{t},j}(S)=(\mathcal{A}, \mathcal{E}_{\mathsf{t},j}, \delta_{\mathsf{t},j})$, where the $(\mathsf{t},j)$-transition function $\delta_{\mathsf{t},j}:\{0,...,q_j-1\}\times \mathcal{A} \to \mathcal{E}_{\mathsf{t},j}$ is defined such that
\begin{equation} \label{Eq:TensorAuto}
\delta_{\mathsf{t},j}(r,a)= b \quad \text{if and only if} \quad [S(a)]\big( \mathsf{t}- D_\mathsf{q}( \mathsf{t} ) + r \mathsf{e}_j \big)=b.
\end{equation}
The $T^j$-patch automaton is the tensor product automaton
\[ \mathcal{G}_{T^j}(S) := \bigotimes_{\mathsf{t}\in T^j}\mathcal{G}_{\mathsf{t},j}.\]
\end{defn}

\begin{figure}[ht]
\centering

\begin{subfigure}[t]{0.45\textwidth}
\centering
\scalebox{0.85}
{
\begin{tikzpicture}
	[
				node distance= 0.1mm, 
				0node/.style={rectangle, draw=black, fill=gray!50!white, thick, minimum size=5mm, font ={$0$} },
				1node/.style={rectangle, draw=black, fill=red!50!white, thick, minimum size=5mm, font ={$1$}},
				2node/.style={rectangle, draw=black, fill=green!50!white, thick, minimum size=5mm, font ={$2$}},
				3node/.style={rectangle, draw=black, fill=blue!50!white, thick, minimum size=5mm, font={$3$}},
				scale=0.9
	]    
	
	\begin{scope}
		\node[0node] (0) {};
	\end{scope}
	
	\begin{scope}[xshift =2.5cm]
		\node[1node] (1) {};
	\end{scope}
	
	\begin{scope}[yshift =-2.5cm]
		\node[2node] (2) {};
	\end{scope}
	
	\begin{scope}[xshift =2.5cm, yshift=-2.5cm]
		\node[3node] (3) {};
	\end{scope}
	
	\draw[-latex, very thick] (2) edge  node[right] {$1$}  (1);	
	\draw[-latex, bend right, very thick, looseness=0.85] (1) edge  node[left] {$1$}  (2);	
	
	\draw[-latex, very thick] (3) edge  node[above] {$0$}  (2);	
	\draw[-latex, bend right, very thick, looseness=0.85] (2) edge  node[below] {$0$}  (3);
			
	\draw[->, very thick] (0) to[out=220, in=140, looseness=5, min distance=0.5cm] node[left] {$0,1$} (0);	
	\draw[->, very thick] (1) to[out=40, in=-40, looseness=5, min distance=0.5cm] node[right] {$0$} (1);
	\draw[->, very thick] (3) to[out=40, in=-40, looseness=5, min distance=0.5cm] node[right] {$1$} (3);

\end{tikzpicture}
}
\caption{The automaton $\mathcal{G}_{ (-1,0),2}(S_{\texttt{table}})$.}
\end{subfigure}
\hfill
\begin{subfigure}[t]{0.45\textwidth}
\centering
\scalebox{0.85}
{
\begin{tikzpicture}
    [
				node distance= 0.1mm, 
				0node/.style={rectangle, draw=black, fill=gray!50!white, thick, minimum size=5mm, font ={$0$} },
				1node/.style={rectangle, draw=black, fill=red!50!white, thick, minimum size=5mm, font ={$1$}},
				2node/.style={rectangle, draw=black, fill=green!50!white, thick, minimum size=5mm, font ={$2$}},
				3node/.style={rectangle, draw=black, fill=blue!50!white, thick, minimum size=5mm, font={$3$}},
				scale=0.9
	]    
	
	\begin{scope}
		\node[0node] (0) {};
	\end{scope}
	
	\begin{scope}[xshift =2.5cm]
		\node[1node] (1) {};
	\end{scope}
	
	\begin{scope}[yshift =-2.5cm]
		\node[2node] (2) {};
	\end{scope}
	
	\begin{scope}[xshift =2.5cm, yshift=-2.5cm]
		\node[3node] (3) {};
	\end{scope}
	
	\draw[-latex, very thick] (3) edge  node[right] {$0$}  (0);	
	\draw[-latex, bend right, very thick, looseness=0.85] (0) edge  node[left] {$0$}  (3);
	
	\draw[-latex, very thick] (0) edge  node[below] {$1$}  (1);	
	\draw[-latex, bend right, very thick, looseness=0.85] (1) edge  node[above] {$1$}  (0);

	\draw[->, very thick] (2) to[out=220, in=140, looseness=5, min distance=0.5cm] node[left] {$0,1$} (2);
	\draw[->, very thick] (1) to[out=40, in=-40, looseness=5, min distance=0.5cm] node[right] {$0$} (1);
	\draw[->, very thick] (3) to[out=40, in=-40, looseness=5, min distance=0.5cm] node[right] {$1$} (3);
\end{tikzpicture}
}
\caption{The automaton $\mathcal{G}_{ (0,0),2}(S_{\texttt{table}})$.}
\end{subfigure}

\caption{The automata comprising $\mathcal{G}_{T^2}(S_{\texttt{table}})$.}
\label{Fig:TableTensor2}
\end{figure}
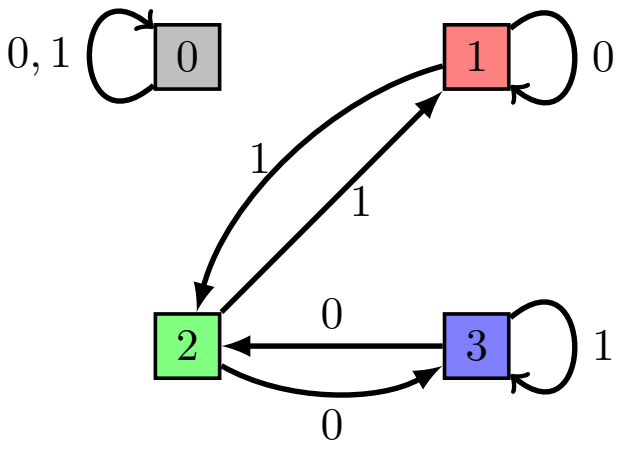
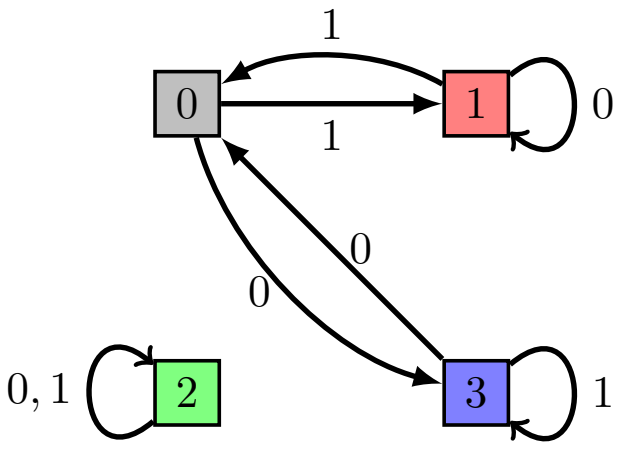

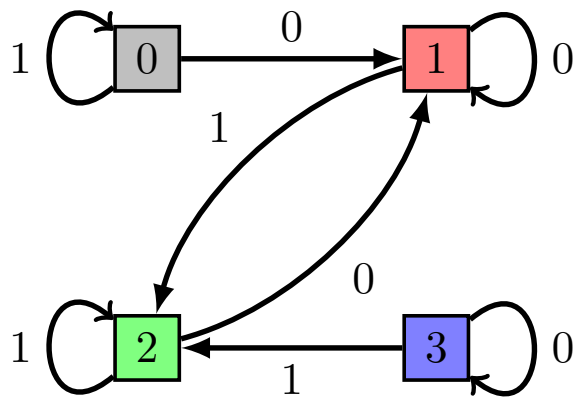
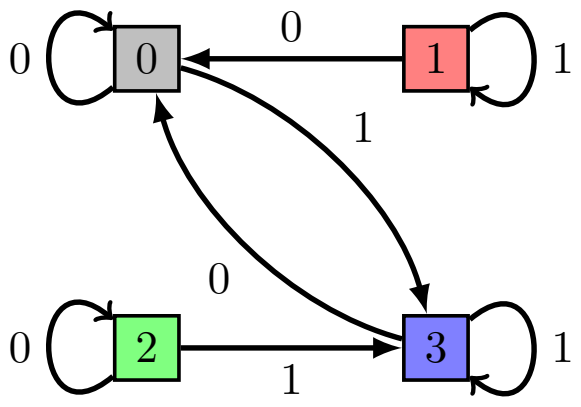
\begin{figure}[ht]
\centering

\begin{subfigure}[t]{0.45\textwidth}
\centering
\scalebox{0.85}
{
\begin{tikzpicture}
	[
				node distance= 0.1mm, 
				0node/.style={rectangle, draw=black, fill=gray!50!white, thick, minimum size=5mm, font ={$0$} },
				1node/.style={rectangle, draw=black, fill=red!50!white, thick, minimum size=5mm, font ={$1$}},
				2node/.style={rectangle, draw=black, fill=green!50!white, thick, minimum size=5mm, font ={$2$}},
				3node/.style={rectangle, draw=black, fill=blue!50!white, thick, minimum size=5mm, font={$3$}},
				scale=0.9
	]    
	
	\begin{scope}
		\node[0node] (0) {};
	\end{scope}
	
	\begin{scope}[xshift =2.5cm]
		\node[1node] (1) {};
	\end{scope}
	
	\begin{scope}[yshift =-2.5cm]
		\node[2node] (2) {};
	\end{scope}
	
	\begin{scope}[xshift =2.5cm, yshift=-2.5cm]
		\node[3node] (3) {};
	\end{scope}
	
	\draw[-latex, bend right, very thick, looseness=0.85] (2) edge  node[below right] {$0$}  (1);	
	\draw[-latex, bend right, very thick, looseness=0.85] (1) edge  node[above left] {$1$}  (2);	
	
	\draw[-latex, very thick] (3) edge  node[below] {$1$}  (2);	
	\draw[-latex, very thick] (0) edge  node[above] {$0$}  (1);

	\draw[->, very thick] (0) to[out=220, in=140, looseness=5, min distance=0.5cm] node[left] {$1$} (0);
	\draw[->, very thick] (2) to[out=220, in=140, looseness=5, min distance=0.5cm] node[left] {$1$} (2);	
	\draw[->, very thick] (1) to[out=40, in=-40, looseness=5, min distance=0.5cm] node[right] {$0$} (1);
	\draw[->, very thick] (3) to[out=40, in=-40, looseness=5, min distance=0.5cm] node[right] {$0$} (3);

\end{tikzpicture}
}
\caption{The automaton $\mathcal{G}_{ (0,-1),1}(S_{\texttt{chair}})$.}
\end{subfigure}
\hfill
\begin{subfigure}[t]{0.45\textwidth}
\centering
\scalebox{0.85}
{
\begin{tikzpicture}
    [
				node distance= 0.1mm, 
				0node/.style={rectangle, draw=black, fill=gray!50!white, thick, minimum size=5mm, font ={$0$} },
				1node/.style={rectangle, draw=black, fill=red!50!white, thick, minimum size=5mm, font ={$1$}},
				2node/.style={rectangle, draw=black, fill=green!50!white, thick, minimum size=5mm, font ={$2$}},
				3node/.style={rectangle, draw=black, fill=blue!50!white, thick, minimum size=5mm, font={$3$}},
				scale=0.9
	]    
	
	\begin{scope}
		\node[0node] (0) {};
	\end{scope}
	
	\begin{scope}[xshift =2.5cm]
		\node[1node] (1) {};
	\end{scope}
	
	\begin{scope}[yshift =-2.5cm]
		\node[2node] (2) {};
	\end{scope}
	
	\begin{scope}[xshift =2.5cm, yshift=-2.5cm]
		\node[3node] (3) {};
	\end{scope}
	
	\draw[-latex, bend left, very thick, looseness=0.85] (0) edge  node[above right] {$1$}  (3);	
	\draw[-latex, bend left, very thick, looseness=0.85] (3) edge  node[below left] {$0$}  (0);
	
	\draw[-latex, very thick] (1) edge  node[above] {$0$}  (0);	
	\draw[-latex, very thick] (2) edge  node[below] {$1$}  (3);

	\draw[->, very thick] (0) to[out=220, in=140, looseness=5, min distance=0.5cm] node[left] {$0$} (0);	
	\draw[->, very thick] (2) to[out=220, in=140, looseness=5, min distance=0.5cm] node[left] {$0$} (2);
	\draw[->, very thick] (1) to[out=40, in=-40, looseness=5, min distance=0.5cm] node[right] {$1$} (1);
	\draw[->, very thick] (3) to[out=40, in=-40, looseness=5, min distance=0.5cm] node[right] {$1$} (3);
\end{tikzpicture}
}
\caption{The automaton $\mathcal{G}_{ (0,0),1}(S_{\texttt{chair}})$.}
\end{subfigure}

\caption{The automata comprising $\mathcal{G}_{T^1}(S_{\texttt{chair}})$.}
\label{Fig:ChairTensor1}
\end{figure}

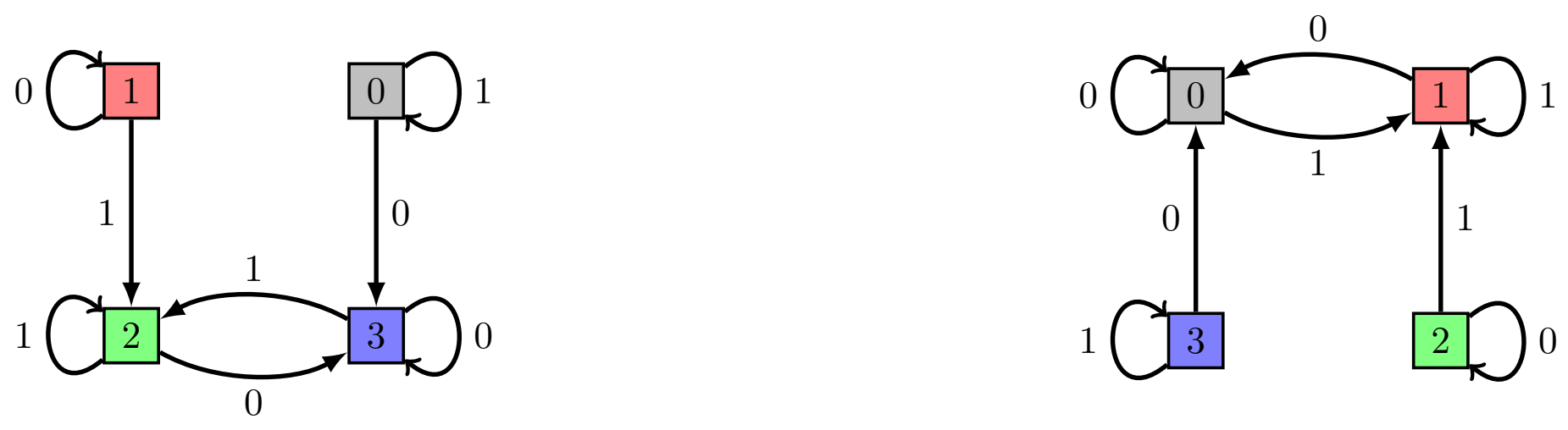
\begin{figure}[ht]
\centering

\begin{subfigure}[t]{0.45\textwidth}
\centering
\scalebox{0.85}
{
\begin{tikzpicture}
	[
				node distance= 0.1mm, 
				0node/.style={rectangle, draw=black, fill=gray!50!white, thick, minimum size=5mm, font ={$0$} },
				1node/.style={rectangle, draw=black, fill=red!50!white, thick, minimum size=5mm, font ={1}},
				2node/.style={rectangle, draw=black, fill=green!50!white, thick, minimum size=5mm, font ={2}},
				3node/.style={rectangle, draw=black, fill=blue!50!white, thick, minimum size=5mm, font={3}},
				scale=0.9
	]    
	
	\begin{scope}
		\node[1node] (1) {};
	\end{scope}
	
	\begin{scope}[xshift =2.5cm]
		\node[0node] (0) {};
	\end{scope}
	
	\begin{scope}[yshift =-2.5cm]
		\node[2node] (2) {};
	\end{scope}
	
	\begin{scope}[xshift =2.5cm, yshift=-2.5cm]
		\node[3node] (3) {};
	\end{scope}
	
	\draw[-latex, bend right, very thick, looseness=0.85] (2) edge  node[below] {$0$}  (3);	
	\draw[-latex, bend right, very thick, looseness=0.85] (3) edge  node[above] {$1$}  (2);	
	
	\draw[-latex, very thick] (1) edge  node[left] {$1$}  (2);	
	\draw[-latex, very thick] (0) edge  node[right] {$0$}  (3);

	\draw[->, very thick] (1) to[out=220, in=140, looseness=5, min distance=0.5cm] node[left] {$0$} (1);
	\draw[->, very thick] (2) to[out=220, in=140, looseness=5, min distance=0.5cm] node[left] {$1$} (2);	
	\draw[->, very thick] (0) to[out=40, in=-40, looseness=5, min distance=0.5cm] node[right] {$1$} (0);
	\draw[->, very thick] (3) to[out=40, in=-40, looseness=5, min distance=0.5cm] node[right] {$0$} (3);

\end{tikzpicture}
}
\caption{The automaton $\mathcal{G}_{ (-1,0),2}(S_{\texttt{chair}})$.}
\end{subfigure}
\hfill
\begin{subfigure}[t]{0.45\textwidth}
\centering
\scalebox{0.85}
{
\begin{tikzpicture}
    [
				node distance= 0.1mm, 
				0node/.style={rectangle, draw=black, fill=gray!50!white, thick, minimum size=5mm, font ={$0$} },
				1node/.style={rectangle, draw=black, fill=red!50!white, thick, minimum size=5mm, font ={$1$}},
				2node/.style={rectangle, draw=black, fill=green!50!white, thick, minimum size=5mm, font ={$2$}},
				3node/.style={rectangle, draw=black, fill=blue!50!white, thick, minimum size=5mm, font={$3$}},
				scale=0.9
	]    
	
	\begin{scope}
		\node[0node] (0) {};
	\end{scope}
	
	\begin{scope}[xshift =2.5cm]
		\node[1node] (1) {};
	\end{scope}
	
	\begin{scope}[yshift =-2.5cm]
		\node[3node] (3) {};
	\end{scope}
	
	\begin{scope}[xshift =2.5cm, yshift=-2.5cm]
		\node[2node] (2) {};
	\end{scope}
	
	\draw[-latex, bend right, very thick, looseness=0.85] (0) edge  node[below] {$1$}  (1);	
	\draw[-latex, bend right, very thick, looseness=0.85] (1) edge  node[above] {$0$}  (0);
	
	\draw[-latex, very thick] (2) edge  node[right] {$1$}  (1);
	\draw[-latex, very thick] (3) edge  node[left] {$0$}  (0);

	\draw[->, very thick] (0) to[out=220, in=140, looseness=5, min distance=0.5cm] node[left] {$0$} (0);	
	\draw[->, very thick] (3) to[out=220, in=140, looseness=5, min distance=0.5cm] node[left] {$1$} (3);
	\draw[->, very thick] (1) to[out=40, in=-40, looseness=5, min distance=0.5cm] node[right] {$1$} (1);
	\draw[->, very thick] (2) to[out=40, in=-40, looseness=5, min distance=0.5cm] node[right] {$0$} (2);
\end{tikzpicture}
}
\caption{The automaton $\mathcal{G}_{ (0,0),2}(S_{\texttt{chair}})$.}
\end{subfigure}

\caption{The automata comprising $\mathcal{G}_{T^2}(S_{\texttt{chair}})$.}
\label{Fig:ChairTensor2}
\end{figure}

For the symbolic table tiling substitution, its $T^1$-patch and $T^2$-patch automata are given in  Figures \ref{Fig:TableTensor1} and \ref{Fig:TableTensor2} respectively. 
Compare the $T^1$-patch and $T^2$-patch automata in Figures~\ref{Fig:TableTensor1} and \ref{Fig:TableTensor2} with the corresponding $1$-st and $2$-nd axis substitution automata in Figures~\ref{Fig:Table-1SubAuto} and \ref{Fig:Table-2SubAuto}, respectively.
For the symbolic chair tiling substitution, its $T^1$-patch and $T^2$-patch automata are given in Figures \ref{Fig:ChairTensor1} and \ref{Fig:ChairTensor2} respectively.

\begin{figure}[ht]
\centering

\begin{subfigure}[t]{0.45\textwidth}
\centering
\scalebox{0.85}
	{
\begin{tikzpicture}
	[
				node distance= 0.1mm, 
				0node/.style={rectangle, draw=black, fill=gray!50!white, thick, minimum size=5mm, font ={$0$} },
				1node/.style={rectangle, draw=black, fill=red!50!white, thick, minimum size=5mm, font ={$1$}},
				2node/.style={rectangle, draw=black, fill=green!50!white, thick, minimum size=5mm, font ={$2$}},
				3node/.style={rectangle, draw=black, fill=blue!50!white, thick, minimum size=5mm, font={$3$}},
				scale=0.9
	]    
	
	\begin{scope}
		\node[0node] (0) {};
	\end{scope}
	
	\begin{scope}[xshift =2.5cm]
		\node[1node] (1) {};
	\end{scope}
	
	\begin{scope}[yshift =-2.5cm]
		\node[2node] (2) {};
	\end{scope}
	
	\begin{scope}[xshift =2.5cm, yshift=-2.5cm]
		\node[3node] (3) {};
	\end{scope}
	
	\draw[-latex, ,bend right, very thick, looseness=0.85] (1) edge  node[above] {$0$}  (0);		
	\draw[-latex, very thick] (0) edge  node[below] {$0$}  (1);	
	\draw[-latex, ,bend right, very thick, looseness=0.85] (2) edge  node[right] {$1$}  (1);	
	\draw[-latex, very thick] (1) edge  node[left] {$1$}  (2);	
			
	\draw[->, very thick] (0) to[out=220, in=140, looseness=5, min distance=0.5cm] node[left] {$1$} (0);
	\draw[->, very thick] (2) to[out=220, in=140, looseness=5, min distance=0.5cm] node[left] {$0$} (2);	
	\draw[->, very thick] (3) to[out=40, in=-40, looseness=5, min distance=0.5cm] node[right] {$0,1$} (3);

\end{tikzpicture}
}
\caption{The automaton $\mathcal{G}_{ (0,-1),1}(S_{\texttt{table}})$.}
\end{subfigure}
\hfill
\begin{subfigure}[t]{0.45\textwidth}
\centering
\scalebox{0.85}
{
\begin{tikzpicture}
    [
				node distance= 0.1mm, 
				0node/.style={rectangle, draw=black, fill=gray!50!white, thick, minimum size=5mm, font ={$0$} },
				1node/.style={rectangle, draw=black, fill=red!50!white, thick, minimum size=5mm, font ={$1$}},
				2node/.style={rectangle, draw=black, fill=green!50!white, thick, minimum size=5mm, font ={$2$}},
				3node/.style={rectangle, draw=black, fill=blue!50!white, thick, minimum size=5mm, font={$3$}},
				scale=0.9
	]    
	
	\begin{scope}
		\node[0node] (0) {};
	\end{scope}
	
	\begin{scope}[xshift =2.5cm]
		\node[1node] (1) {};
	\end{scope}
	
	\begin{scope}[yshift =-2.5cm]
		\node[2node] (2) {};
	\end{scope}
	
	\begin{scope}[xshift =2.5cm, yshift=-2.5cm]
		\node[3node] (3) {};
	\end{scope}
	
	\draw[-latex ,bend right, very thick, looseness=0.85] (3) edge  node[right] {$0$}  (0);	
	\draw[-latex, very thick] (0) edge  node[left] {$0$}  (3);		
	\draw[-latex, very thick] (3) edge  node[above] {$1$}  (2);	
	\draw[-latex, bend right, very thick, looseness=1.25] (2) edge  node[below] {$1$}  (3);	
	
	\draw[->, very thick] (0) to[out=220, in=140, looseness=5, min distance=0.5cm] node[left] {$1$} (0);
	\draw[->, very thick] (2) to[out=220, in=140, looseness=5, min distance=0.5cm] node[left] {$0$} (2);
	\draw[->, very thick] (1) to[out=40, in=-40, looseness=5, min distance=0.5cm] node[right] {$0,1$} (1);
\end{tikzpicture}
}
\caption{The automaton $\mathcal{G}_{ (0,0),1}(S_{\texttt{table}})$.}
\end{subfigure}

\caption{The automata comprising $\mathcal{G}_{T^1}(S_{\texttt{table}})$.}
\label{Fig:TableTensor1}
\end{figure}
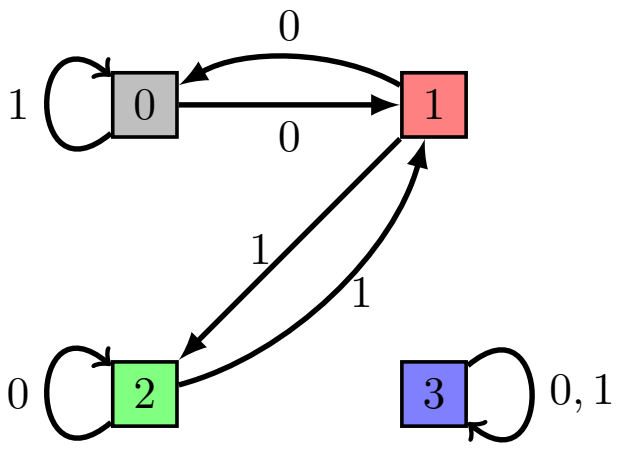
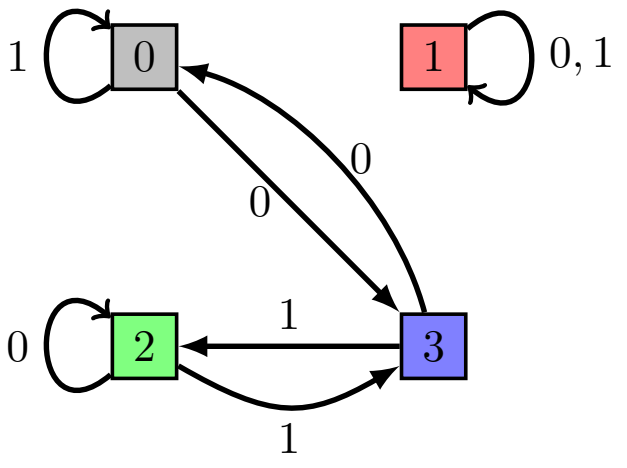

The correspondence between the $T^j$-patch automaton and the $j$-th axis substitution automaton holds in complete generality.  
The vertex set of $\mathcal{G}_{T^j}(S)$ is canonically identified with $\mathcal{A}^{T^j}$ by $(a_\mathsf{t})_{\mathsf{t}\in T^j}\mapsto P$ with $P(\mathsf{t})=a_{\mathsf{t}}$. 
\begin{lemma} \label{Lem:TensorCong}
 Under the canonical identification $\prod_{\mathsf{t}\in T^j } \mathcal{A} \cong \mathcal{A}^{T^j}$, there is a directed edge from $P$ to $P'$ in $\mathcal{G}_{T^j}(S)$ labelled by  $0\leq r\leq q_j-1$  if and only if 
 \begin{equation} \label{Eq:TjIdentity}
 	S(P)|_{r\mathsf{e}_j+T^j}=P'. 
 \end{equation}
 In particular, the $j$-th axis substitution automaton from Definition \ref{Defn:SubAutoJ} is isomorphic to the tensor product automaton $\bigotimes_{\mathsf{t}\in T^j} \mathcal{G}_{\mathsf{t},j}(S)$ from Definition \ref{Defn:TensorAuto}.
\end{lemma}

\begin{proof}
	For every $\mathsf{t}\in T^j$, recall that the transition function $\delta_{\mathsf t,j}$ is defined by \eqref{Eq:TensorAuto}.
	Let $P\cong (a_{\mathsf{t}})_{\mathsf{t}\in T^j}$ and $P'\cong (b_{\mathsf{t}})_{\mathsf{t}\in T^j}$.
	Then there is an edge  from $P$ to $P'$ labelled by $r\in \{0,...,q_j-1\}$	in $\mathcal{G}_{T^j}(S)$ if and only if
	\[
	\delta_{\mathsf t,j}(r,a_{\mathsf t})=
	b_{\mathsf t} \qquad\text{for all }\mathsf t\in T^j.
	\]
	By the definition of \(\delta_{\mathsf t,j}\), this is 	equivalent to
	\[
	[S(P)](\mathsf t+r\mathsf e_j)=P'(\mathsf t)
	\qquad\text{for all }\mathsf t\in T^j,
	\]
	which is precisely \eqref{Eq:TjIdentity}.
\end{proof}

\begin{proof}[Proof of Proposition \ref{Prop:PureLineDefect}]
	Following Corollary~\ref{Cor:PL-Defect}, it suffices to show that the $j$-th axis defect automaton is a complete automaton. Hence, it suffices to show that for every illegal $T^j$-patch $P$ and every label $r\in\{0,\ldots,q_j-1\}$,$$\delta_j(r,P)\notin W(S).$$
	Since $S$ is bijective along the $j$-th axis, for every 
	$\mathsf{r}\in K_{\mathsf{q}}\cap \{ (x_1,...,x_d) :x_j \in \{0,q_j-1\} \}$ and every $a\in \mathcal{A}$, there exists a unique letter $b=b(\mathsf{r},a)\in \mathcal{A}$ such that $[S_0(a)](\mathsf{r})=b$. Moreover, every  $\mathsf{r}\in K_{\mathsf{q}}\cap \{ (x_1,...,x_d) :x_j \in \{0,q_j-1\} \}$ can be written uniquely as 
	$$\mathsf{r}= \mathsf{t} - D_{\mathsf{q}}(t)+r \mathsf{e_j} \quad \text{where} \quad  \mathsf{t}\in T^j\; \text{and} \; r\in \{ 0,...,q_j-1 \}.$$
	Hence,	$S$ is bijective along the $j$-th axis if and only if
	\[ |\{ b\in \mathcal{A}: \delta_{\mathsf{t},j}(a,r)=b \}|=1 \quad \text{for every} \quad \mathsf{t}\in T^j, \; a\in \mathcal{A} \;\; \text{and} \; \; r\in\{0,...,q_j-1\}. \]
	A simple induction argument on $k\in \mathbb{N}$, shows that
	$S$ is bijective along the $j$-th axis if and only if
	\[ |\{ b\in \mathcal{A}: \delta_{\mathsf{t},j}(a,r_1...r_k)=b \}|=1 \quad \text{for every} \quad \mathsf{t}\in T^j, \; a\in \mathcal{A} \;\; \text{and} \; \; r_1,...,r_k\in\{0,...,q_j-1\}. \]
	Fix $r\in \{0,...,q_j-1\}$ and $\mathsf{t}\in T^j$. For every $a\in \mathcal{A}$, the sequence 
	$$\{\delta_{\mathsf{t,j}}\big(a,r^n\big)\}_{n=0}^\infty \quad \text{is eventually periodic,}$$ since $\mathcal{A}$ is finite and the transitions are deterministic. Moreover, bijectivity implies that each transition admits a unique predecessor, so the sequence is in fact periodic. 
	It follows from Lemma \ref{Lem:TensorCong} that the sequence $\delta_j\big(P,r^n\big)$ is also periodic. Therfore, the definition of the $j$-th axis susbtitution automaton implies that $P\prec S^{N}( \delta_j(r,P) )$ for some $N\in \mathbb{N}$. Since $P\notin W(S)$ and $W(S)$ is downwards closed under the subpatch relation, it follows that $\delta_j(r,P)\notin W(S)$, completing the proof.
\end{proof}

\subsection{Powers of automata and powers of substitutions} \label{Sect:SubstPower}

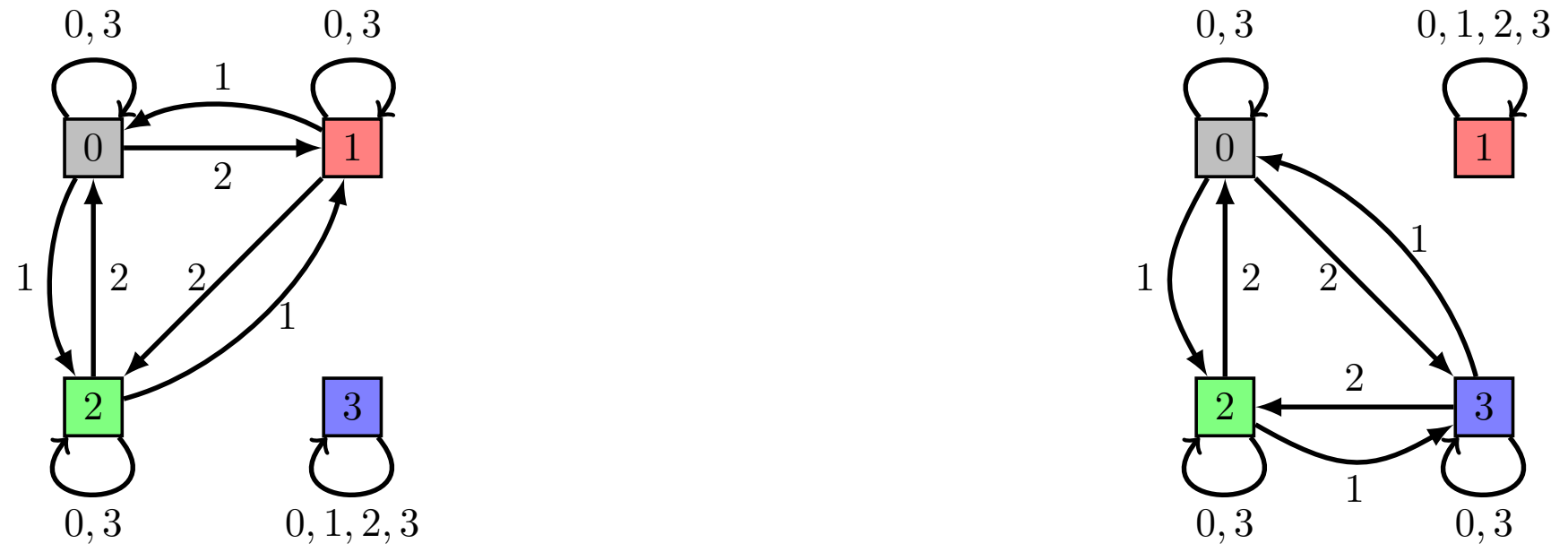
\begin{figure}[ht]
\centering

\begin{subfigure}[t]{0.45\textwidth}
\centering
\scalebox{0.85}
{
\begin{tikzpicture}
	[
				node distance= 0.1mm, 
				0node/.style={rectangle, draw=black, fill=gray!50!white, thick, minimum size=5mm, font ={$0$} },
				1node/.style={rectangle, draw=black, fill=red!50!white, thick, minimum size=5mm, font ={$1$}},
				2node/.style={rectangle, draw=black, fill=green!50!white, thick, minimum size=5mm, font ={$2$}},
				3node/.style={rectangle, draw=black, fill=blue!50!white, thick, minimum size=5mm, font={$3$}},
				scale=0.9
	]    
	
	\begin{scope}
		\node[0node] (0) {};
	\end{scope}
	
	\begin{scope}[xshift =2.5cm]
		\node[1node] (1) {};
	\end{scope}
	
	\begin{scope}[yshift =-2.5cm]
		\node[2node] (2) {};
	\end{scope}
	
	\begin{scope}[xshift =2.5cm, yshift=-2.5cm]
		\node[3node] (3) {};
	\end{scope}
	
	\draw[-latex, ,bend right, very thick, looseness=0.85] (1) edge  node[above] {$1$}  (0);		
	\draw[-latex, very thick] (0) edge  node[below] {$2$}  (1);	
	\draw[-latex, ,bend right, very thick, looseness=0.85] (2) edge  node[right] {$1$}  (1);	
	\draw[-latex, very thick] (1) edge  node[left] {$2$}  (2);
	\draw[-latex, ,bend right, very thick, looseness=0.85] (0) edge  node[left] {$1$}  (2);		
	\draw[-latex, very thick] (2) edge  node[right] {$2$}  (0);
			
	\draw[->, very thick] (0) to[out=130, in=50, looseness=5, min distance=0.5cm] node[above] {$0,3$} (0);
	\draw[->, very thick] (1) to[out=130, in=50, looseness=5, min distance=0.5cm] node[above] {$0,3$} (1);
	\draw[->, very thick] (2) to[out=310, in=230, looseness=5, min distance=0.5cm] node[below] {$0,3$} (2);	
	\draw[->, very thick] (3) to[out=310, in=230, looseness=5, min distance=0.5cm] node[below] {$0,1,2,3$} (3);

\end{tikzpicture}
}
\caption{The automaton $\mathcal{G}_{ (0,-1),1}(S^2_{\texttt{table}})$.}
\end{subfigure}
\hfill
\begin{subfigure}[t]{0.45\textwidth}
\centering
\scalebox{0.85}
{
\begin{tikzpicture}
    [
				node distance= 0.1mm, 
				0node/.style={rectangle, draw=black, fill=gray!50!white, thick, minimum size=5mm, font ={$0$} },
				1node/.style={rectangle, draw=black, fill=red!50!white, thick, minimum size=5mm, font ={$1$}},
				2node/.style={rectangle, draw=black, fill=green!50!white, thick, minimum size=5mm, font ={$2$}},
				3node/.style={rectangle, draw=black, fill=blue!50!white, thick, minimum size=5mm, font={$3$}},
				scale=0.9
	]    
	
	\begin{scope}
		\node[0node] (0) {};
	\end{scope}
	
	\begin{scope}[xshift =2.5cm]
		\node[1node] (1) {};
	\end{scope}
	
	\begin{scope}[yshift =-2.5cm]
		\node[2node] (2) {};
	\end{scope}
	
	\begin{scope}[xshift =2.5cm, yshift=-2.5cm]
		\node[3node] (3) {};
	\end{scope}
	
	\draw[-latex ,bend right, very thick, looseness=0.85] (3) edge  node[right] {$1$}  (0);	
	\draw[-latex, very thick] (0) edge  node[left] {$2$}  (3);		
	\draw[-latex, very thick] (3) edge  node[above] {$2$}  (2);	
	\draw[-latex, bend right, very thick, looseness=1.25] (2) edge  node[below] {$1$}  (3);	
	\draw[-latex, very thick] (2) edge  node[right] {$2$}  (0);
	\draw[-latex, bend right, very thick, looseness=1.25] (0) edge  node[left] {$1$}  (2);	
	
	\draw[->, very thick] (0) to[out=130, in=50, looseness=5, min distance=0.5cm] node[above] {$0,3$} (0);
	\draw[->, very thick] (1) to[out=130, in=50, looseness=5, min distance=0.5cm] node[above] {$0,1,2,3$} (1);
	\draw[->, very thick] (2) to[out=310, in=230, looseness=5, min distance=0.5cm] node[below] {$0,3$} (2);	
	\draw[->, very thick] (3) to[out=310, in=230, looseness=5, min distance=0.5cm] node[below] {$0,3$} (3);
\end{tikzpicture}
}
\caption{The automaton $\mathcal{G}_{ (0,0),1}(S^2_{\texttt{table}})$.}
\end{subfigure}

\caption{The automata comprising $\mathcal{G}_{T^1}(S^2_{\texttt{table}})$.}
\label{Fig:TableTensor1Square}
\end{figure}

\begin{figure}[ht]
\centering

\begin{subfigure}[t]{0.45\textwidth}
\centering
\scalebox{0.85}
{
\begin{tikzpicture}
	[
				node distance= 0.1mm, 
				0node/.style={rectangle, draw=black, fill=gray!50!white, thick, minimum size=5mm, font ={$0$} },
				1node/.style={rectangle, draw=black, fill=red!50!white, thick, minimum size=5mm, font ={$1$}},
				2node/.style={rectangle, draw=black, fill=green!50!white, thick, minimum size=5mm, font ={$2$}},
				3node/.style={rectangle, draw=black, fill=blue!50!white, thick, minimum size=5mm, font={$3$}},
				scale=0.9
	]    
	
	\begin{scope}
		\node[0node] (0) {};
	\end{scope}
	
	\begin{scope}[xshift =2.5cm]
		\node[1node] (1) {};
	\end{scope}
	
	\begin{scope}[yshift =-2.5cm]
		\node[2node] (2) {};
	\end{scope}
	
	\begin{scope}[xshift =2.5cm, yshift=-2.5cm]
		\node[3node] (3) {};
	\end{scope}
	
	\draw[-latex, ,bend right, very thick, looseness=0.85] (3) edge  node[right] {$1$}  (1);		
	\draw[-latex, very thick] (1) edge  node[left] {$2$}  (3);	
	\draw[-latex, ,bend right, very thick, looseness=0.85] (2) edge  node[below] {$1$}  (3);	
	\draw[-latex, very thick] (3) edge  node[above] {$2$}  (2);
	
	\draw[-latex, very thick] (2) edge  node[below] {$2$}  (1);
	\draw[-latex, bend right, very thick, looseness=1.25] (1) edge  node[left] {$1$}  (2);	
			
	\draw[->, very thick] (0) to[out=130, in=50, looseness=5, min distance=0.5cm] node[above] {$0,...,3$} (0);
	\draw[->, very thick] (1) to[out=130, in=50, looseness=5, min distance=0.5cm] node[above] {$0,3$} (1);
	\draw[->, very thick] (2) to[out=310, in=230, looseness=5, min distance=0.5cm] node[below] {$0,3$} (2);	
	\draw[->, very thick] (3) to[out=310, in=230, looseness=5, min distance=0.5cm] node[below] {$0,3$} (3);

\end{tikzpicture}
}
\caption{The automaton $\mathcal{G}_{ (-1,0),2}(S^2_{\texttt{table}})$.}
\end{subfigure}
\hfill
\begin{subfigure}[t]{0.45\textwidth}
\centering
\scalebox{0.85}
{
\begin{tikzpicture}
    [
				node distance= 0.1mm, 
				0node/.style={rectangle, draw=black, fill=gray!50!white, thick, minimum size=5mm, font ={$0$} },
				1node/.style={rectangle, draw=black, fill=red!50!white, thick, minimum size=5mm, font ={$1$}},
				2node/.style={rectangle, draw=black, fill=green!50!white, thick, minimum size=5mm, font ={$2$}},
				3node/.style={rectangle, draw=black, fill=blue!50!white, thick, minimum size=5mm, font={$3$}},
				scale=0.9
	]    
	
	\begin{scope}
		\node[0node] (0) {};
	\end{scope}
	
	\begin{scope}[xshift =2.5cm]
		\node[1node] (1) {};
	\end{scope}
	
	\begin{scope}[yshift =-2.5cm]
		\node[2node] (2) {};
	\end{scope}
	
	\begin{scope}[xshift =2.5cm, yshift=-2.5cm]
		\node[3node] (3) {};
	\end{scope}
	
	\draw[-latex ,bend right, very thick, looseness=0.85] (3) edge  node[right] {$1$}  (1);	
	\draw[-latex, very thick] (1) edge  node[left] {$2$}  (3);		
	\draw[-latex, very thick] (0) edge  node[below] {$2$}  (1);	
	\draw[-latex, bend right, very thick, looseness=1.25] (1) edge  node[above] {$1$}  (0);	
	
	\draw[-latex, very thick] (3) edge  node[right] {$2$}  (0);
	\draw[-latex, bend right, very thick, looseness=1.25] (0) edge  node[left] {$1$}  (3);
	
	\draw[->, very thick] (0) to[out=130, in=50, looseness=5, min distance=0.5cm] node[above] {$0,3$} (0);
	\draw[->, very thick] (1) to[out=130, in=50, looseness=5, min distance=0.5cm] node[above] {$0,3$} (1);
	\draw[->, very thick] (2) to[out=310, in=230, looseness=5, min distance=0.5cm] node[below] {$0,...,3$} (2);	
	\draw[->, very thick] (3) to[out=310, in=230, looseness=5, min distance=0.5cm] node[below] {$0,3$} (3);
\end{tikzpicture}
}
\caption{The automaton $\mathcal{G}_{ (0,0),2}(S^2_{\texttt{table}})$.}
\end{subfigure}

\caption{The automata comprising $\mathcal{G}_{T^2}(S^2_{\texttt{table}})$.}
\label{Fig:TableTensor2Square}
\end{figure}
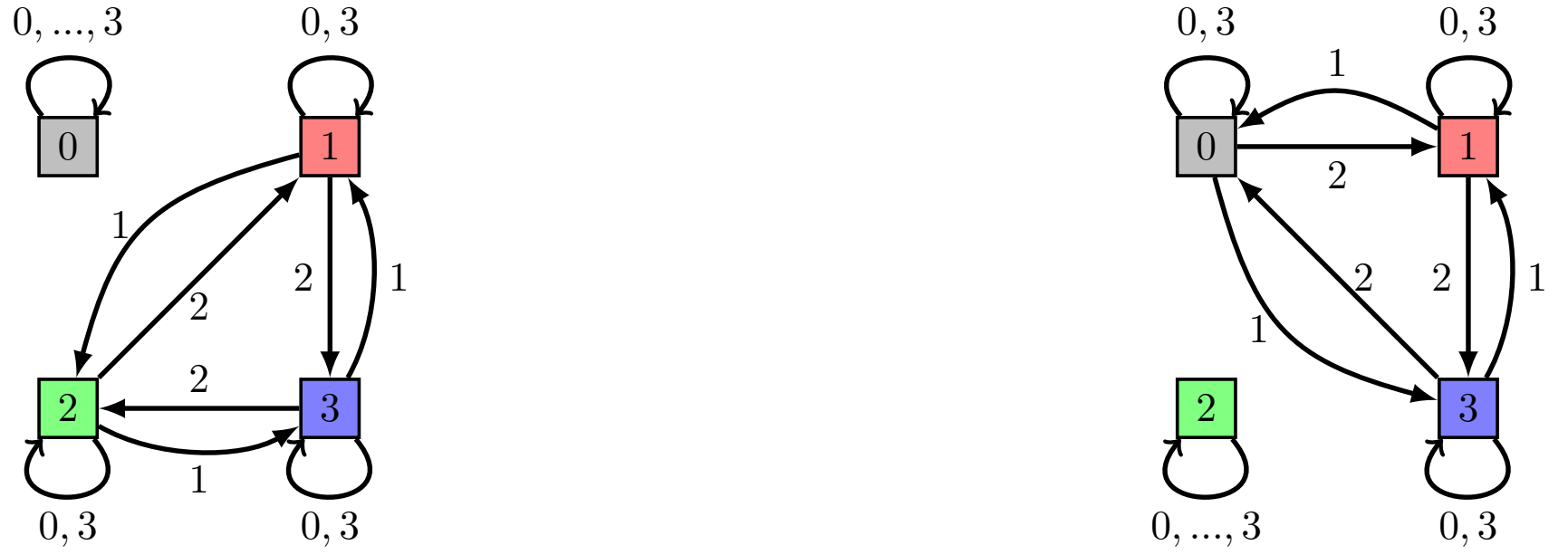

It is often useful to consider powers of a substitution. Given a block substitution $S$, its $k$-th power of a substitution $S$ is simply $S^k$. 
In particular, one may study configurations $\rho\in \mathcal{A}^{\mathbb{Z}^d}$ satisfying $S^k(\rho)=\rho$. 
We therefore describe how the automata introduced in this paper can be constructed for the substitution $S^k$. As a motivating example, the symbolic table substitution admits several fixed points , the symbolic table tiling substitution admits several fixed points $\rho\in \{0,...,3\}^{\mathbb{Z}^2}$ satisfying $S^2_{\texttt{table}}(\rho)=\rho$ but $S_{\texttt{table}}(\rho)\neq \rho$. 
Consider the  automata associated with $S^2_{\texttt{table}}$ introduced in Definition \ref{Defn:TensorAuto}, shown in Figures \ref{Fig:TableTensor1Square} and \ref{Fig:TableTensor2Square}. We explain how these automata can be constructed directly from the axis substitution automata of $S_{\texttt{table}}$.

For a directed graph $(\mathcal{V}, \mathcal{E})$, its $k$-th power is the directed graph $(\mathcal{V}, \mathcal{E}^{k})$, where $(v,u)\in \mathcal{E}^{k}$ if and only if there exists a directed walk of length $k$ starting at $v$ and ending at $u$. Thus, $\mathcal{E}^k$ records reachability in exactly $k$ steps. 
For a determinstic automaton $\mathcal{G}=(\mathcal{V}, \mathcal{E}, \delta)$, its $k$-th power is the  deterministic automaton $\mathcal{G}^{(k)}:=(\mathcal{V}, \mathcal{E}^k, \delta^{(k)})$, where the transition partial function $\delta^{(k)}:L^k \times \mathcal{V}\to \mathcal{V}$ is given by
$$ \delta^{(k)}\big( \ell_1...\ell_k, v \big):=  \delta(\ell_1...\ell_k,v)$$
for all $v\in \mathcal{V}$ and $\ell_1,...,\ell_k\in L^k$ whenever the right hand side is well defined. 
We show that the $j$-th axis substitution and defect automata of $S^k$ coincide with the $k$-th powers of the appropriate automata for $S$.
Note that if $S$ is induced from a substitution rule of length $\mathsf{q}=(q_1,...,q_d)$, then $S^k$ is induced from a substitution rule of length $\mathsf{q}^k:=(q_1^k,...,q_d^k)$.  For every $0\leq \mathbf{r}^{(k)} \leq q_j^k-1$, there exist unique $r^{(1)},...,r^{(k)}\in \{0,..., q_j-1\}$ such that
\begin{equation} \label{Eq:LabelPower}
\mathbf{r}^{(k)}= \sum_{i=1}^k r^{(i)} q_j^{k-i}. 
\end{equation}
Hence, \eqref{Eq:LabelPower} yields a natural identification
\[
\{0,\ldots,q_j^k-1\}
\cong
\{0,\ldots,q_j-1\}^k,
\]
namely the correspondence between an integer and its base-$q_j$ expansion of length $k$.
\begin{lemma} \label{Lem:PowerCong}
	Let $S$ be a block substitution of length $\mathsf{q}=(q_1,...,q_d)$ and fix $1\leq j \leq d$. Then the following assertions hold. 
	\begin{enumerate}[(a)]
	\item For every $\mathsf{t}\in T^j$, $a,b\in \mathcal{A}$ and $r^{(1)}...r^{(k)}\in \{0,...,q_j-1 \}$, 
	\[  \delta^{(k)}_{\mathsf{t},j}(r^{(1)}...r^{(k)},a)=b \quad \text{if and only if} \quad [S^k(a)]\big( \mathsf{t}- D_\mathsf{q}^k(\mathsf{t})+\mathbf{r}^{(k)}\mathsf{e}_j\big)=b,\]
	where $\mathbf{r}^{(k)}$ is given by \eqref{Eq:LabelPower}. In particular, the $k$-th power automaton of $\mathcal{G}_{\mathsf{t},j}(S)$ is isomorphic to $\mathcal{G}_{\mathsf{t},j}(S^k)$.
	\item For every $P,Q\in \mathcal{A}^{T^j}$ and $r^{(1)}...r^{(k)}\in \{0,...,q_j-1 \}$,
	\[  \delta^{(k)}_{j}(r^{(1)}...r^{(k)},P)=Q \quad \text{if and only if} \quad \big(S^k(P)\big)|_{\mathbf{r}^{(k)}\mathsf{e}_j+T^j}=Q,\]
	where $\mathbf{r}^{(k)}$ is given by \eqref{Eq:LabelPower}.  In particular, the $k$-th power automaton of $\mathcal{G}_{j}(S)$ is isomorphic to $\mathcal{G}_{j}(S^k)$.
	\end{enumerate}	
\end{lemma}

\begin{example} \label{Examp:DefectForm2}
	We illustrate how to utilize the mechanisms introduced in Section \ref{Sect:Simpl} to compute a structural defect in an alternative way, revisiting Example~\ref{Examp:DefectForm}. 
	We again consider the symbolic table-tiling substitution, the patch  $P_1$ shown in Figure \ref{Fig:P1Patch}, and the fixed point $\rho_{P_1}$ satisfying $S^2_{\mathtt{table}}(\rho_{P_1})=\rho_{P_1}$. Our goal is to compute $P_1':=\rho_{P_1}|_{m'\mathsf{e}_1+T^1}$ for $m'=6$. A short computation shows that the $4$-coding of $m'=6$ is $(0;1,2)$. Furthermore, since $S^2_{\mathtt{table}}(\rho_{P_1})=\rho_{P_1}$, the modified address has $\hat{r}=0$ delay substitution steps. 
	Using the identification in Lemma \ref{Lem:PowerCong}~(a), we employ the automaton $\mathcal{G}_{(0,0),1}(S^2_{\texttt{table}})$ in Figure~\ref{Fig:TableTensor1Square}~(B) to determine $\rho_{P_1}(m',0)$. The corresponding derivation is presented in Figure~\ref{Fig:m'DerivGenU}. Similarly, in Figure~\ref{Fig:m'DerivGenB} we determine $\rho_{P_1}(m',-1)$, using the automaton $\mathcal{G}_{(0,-1),1}(S^2_{\texttt{table}})$ in Figure~\ref{Fig:TableTensor1Square}~(A).
	
\begin{multicols}{2}

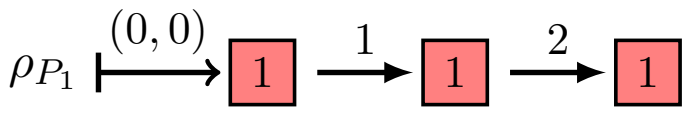
\begin{figure}[H]
	\centering
	\scalebox{0.85}
	{
	\begin{tikzpicture}
		[
				node distance= 0.1mm, 
				0node/.style={rectangle, draw=black, fill=gray!50!white, thick, minimum size=5mm, font ={$0$} },
				1node/.style={rectangle, draw=black, fill=red!50!white, thick, minimum size=5mm, font ={$1$}},
				2node/.style={rectangle, draw=black, fill=green!50!white, thick, minimum size=5mm, font ={$2$}},
				3node/.style={rectangle, draw=black, fill=blue!50!white, thick, minimum size=5mm, font={$3$}},
				]
				
			\begin{scope}[xshift=-2cm]
			\draw[draw=none,xshift=0.3cm, yshift=-0.2cm] (0,0) circle  (0.15cm)node (R-1) {\large $\rho_{P_1}$};
			\end{scope}				
				
			\begin{scope}[xshift=0cm]
			\draw[draw=none,xshift=-0.2cm, yshift=-0.2cm] (0,0) circle  (0.15cm)node (L0) {};
			\draw[draw=none,xshift=0.3cm, yshift=-0.2cm] (0,0) circle  (0.15cm)node (R0) {};
				\node[1node, yshift= -0.2cm] (U0) {};
			\end{scope}
			
			\begin{scope}[xshift=1.5cm]
			\draw[draw=none,xshift=-0.2cm, yshift=-0.2cm] (0,0) circle  (0.15cm)node (L1) {};
			\draw[draw=none,xshift=0.3cm, yshift=-0.2cm] (0,0) circle  (0.15cm)node (R1) {};
				\node[1node, yshift= -0.2cm] (U1) {};
			\end{scope}	
			
			\begin{scope}[xshift=3cm]
			\draw[draw=none,xshift=-0.2cm, yshift=-0.2cm] (0,0) circle  (0.15cm)node (L2) {};
			\draw[draw=none,xshift=0.3cm, yshift=-0.2cm] (0,0) circle  (0.15cm)node (R2) {};
				\node[1node, yshift= -0.2cm] (U2) {};
			\end{scope}

			\draw[|->, very thick] (R-1) edge node[above] {$(0,0)$} (L0);		
			\draw[-latex, very thick] (R0) edge node[above] {$1$} (L1);				
			\draw[-latex, very thick] (R1) edge node[above] {$2$} (L2);
		
	\end{tikzpicture}
	}
	\caption{\centering The derivation for $\rho_{P_1}(m',0)$.}
	\label{Fig:m'DerivGenU}
\end{figure}

\columnbreak

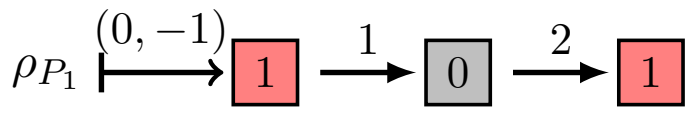
\begin{figure}[H]
	\centering
	\scalebox{0.85}
	{
	\begin{tikzpicture}
		[
				node distance= 0.1mm, 
				0node/.style={rectangle, draw=black, fill=gray!50!white, thick, minimum size=5mm, font ={$0$} },
				1node/.style={rectangle, draw=black, fill=red!50!white, thick, minimum size=5mm, font ={$1$}},
				2node/.style={rectangle, draw=black, fill=green!50!white, thick, minimum size=5mm, font ={$2$}},
				3node/.style={rectangle, draw=black, fill=blue!50!white, thick, minimum size=5mm, font={$3$}},
				]
				
			\begin{scope}[xshift=-2cm]
			\draw[draw=none,xshift=0.3cm, yshift=-0.2cm] (0,0) circle  (0.15cm)node (R-1) {\large $\rho_{P_1}$};
			\end{scope}				
				
			\begin{scope}[xshift=0cm]
			\draw[draw=none,xshift=-0.2cm, yshift=-0.2cm] (0,0) circle  (0.15cm)node (L0) {};
			\draw[draw=none,xshift=0.3cm, yshift=-0.2cm] (0,0) circle  (0.15cm)node (R0) {};
				\node[1node, yshift= -0.2cm] (U0) {};
			\end{scope}
			
			\begin{scope}[xshift=1.5cm]
			\draw[draw=none,xshift=-0.2cm, yshift=-0.2cm] (0,0) circle  (0.15cm)node (L1) {};
			\draw[draw=none,xshift=0.3cm, yshift=-0.2cm] (0,0) circle  (0.15cm)node (R1) {};
				\node[0node, yshift= -0.2cm] (U1) {};
			\end{scope}	
			
			\begin{scope}[xshift=3cm]
			\draw[draw=none,xshift=-0.2cm, yshift=-0.2cm] (0,0) circle  (0.15cm)node (L2) {};
			\draw[draw=none,xshift=0.3cm, yshift=-0.2cm] (0,0) circle  (0.15cm)node (R2) {};
				\node[1node, yshift= -0.2cm] (U2) {};
			\end{scope}

			\draw[|->, very thick] (R-1) edge node[above] {$(0,-1)$} (L0);		
			\draw[-latex, very thick] (R0) edge node[above] {$1$} (L1);				
			\draw[-latex, very thick] (R1) edge node[above] {$2$} (L2);
		
	\end{tikzpicture}
	}
	\caption{\centering The derivation for $\rho_{P_1}(m',-1)$.}
	\label{Fig:m'DerivGenB}
\end{figure}

\end{multicols}

The resulting patch coincides with the one obtained in  Example~\ref{Examp:DefectForm}.

\end{example}

\begin{proof}[Proof of Lemma \ref{Lem:PowerCong}]$~$\\
(a) We prove the claim by induction on $k\in \mathbb{N}$. 
For every $\mathsf{t}\in T^j$, recall that the transition function $\delta^{(k)}_{\mathsf{t},j}$ satisfies \eqref{Eq:TensorAuto}. Hence, the claim follows directly from \eqref{Eq:TensorAuto} when  $k=1$. 
Assume now that the claim holds for some $k\in \mathbb{N}$. 
Then, by \eqref{Eq:TensorAuto} and the induction hypothesis,
\[ \begin{aligned} [S^{k+1}(a)]\big( \mathsf{t}- D_{\mathsf{q}}^{k+1}(\mathsf{t})+ \mathbf{r}^{(k+1)} \mathsf{e}_j \big)&= 
\Big[ S\Big([S^k(a)] \big( \mathsf{t}- D_{\mathsf{q}}(\mathsf{t})+ \mathbf{r}^{(k)}\mathsf{e}_j \big) \Big) \Big]\big(\mathsf{t}- D_{\mathsf{q}}(\mathsf{t})+ r^{(k+1)}\mathsf{e}_j \big) \\ &= \Big[ S\Big( \delta_{\mathsf{t},j}(r^{(1)}...r^{(k)},a) \Big) \Big] \big(\mathsf{t}- D_{\mathsf{q}}(\mathsf{t})+ r^{(k+1)}\mathsf{e}_j \big)\\ 
&= \delta_{\mathsf{t},j} \Big( r^{(k+1)}, \delta_{\mathsf{t},j}(r^{(1)}...r^{(k)},a) \Big) \\
&= \delta^{(k+1)}_{\mathsf{t},j}(r^{(1)}...r^{(k)}r^{(k+1)},a), \end{aligned}  \]
which completes the induction.
\newline
(b) The equality holds if and only if $(P,Q)$ is an edge in $\mathcal{G}_j(S^k)$ labelled by $\mathbf{r}^{(k)}$. By Lemma \ref{Lem:TensorCong}, we know that the automata $\bigotimes_{\mathsf{t}\in T^j} \mathcal{G}_{\mathsf{t},j}(S^k)$  and $\mathcal{G}_j(S^k)$ are isomorphic. Hence, the desired equivalence follows from the previous part.

\end{proof}
\subsection*{Acknowledgements}
I wish to sincerely thank  Ram Band, Siegfried Beckus and Felix Pogorzelski. Their collaboration and discussions with me on \cite{BBPT24, BBPT26} led to this paper. I am also grateful to Philipp Bartmann and Alan Lew for discussions which contributed to this work.

This research was supported by the Israel Science Foundation, grant No. 1647/23 (PI: B. Solomyak).

\bibliographystyle{alpha}
\bibliography{AddressVer3}

\end{document}